\documentclass[12pt,oneside]{amsart}
\usepackage[T1]{fontenc}
\usepackage[ruled,vlined]{algorithm2e}
\usepackage{mathtools,bm,amsthm,amsmath,amsfonts,amssymb,color,epic}
\usepackage{multirow,comment}
\usepackage[export]{adjustbox}
\usepackage{tikz,pgfplots,nicefrac}
\pgfplotsset{compat=newest}
\usetikzlibrary{arrows,snakes,backgrounds,positioning,calc}
\usepackage[left=3cm,right=3cm,top=2.5cm,bottom=2.5cm]{geometry}
\theoremstyle{plain}

\usepackage{subfig}
\theoremstyle{plain}             
\newtheorem{theorem}{Theorem}[section]
\newtheorem{lemma}[theorem]{Lemma}

\def\Letters{A,B,C,D,E,F,G,H,I,J,K,L,M,N,O,P,Q,R,S,T,U,V,W,X,Y,Z}
\makeatletter
\@for \@l:=\Letters \do{%
  \expandafter\edef\csname\@l bb\endcsname{\noexpand\ensuremath{\noexpand\mathbb{\@l}}}%
  \expandafter\edef\csname\@l bf\endcsname{{\noexpand\bs \@l}}%
  \expandafter\edef\csname\@l cal\endcsname{\noexpand\ensuremath{\noexpand\mathcal{\@l}}}%
  \expandafter\edef\csname\@l eu\endcsname{\noexpand\ensuremath{\noexpand\EuScript{\@l}}}%
  \expandafter\edef\csname\@l frak\endcsname{\noexpand\ensuremath{\noexpand\mathfrak{\@l}}}%
  \expandafter\edef\csname\@l rm\endcsname{{\noexpand\rm \@l}}%
  \expandafter\edef\csname\@l scr\endcsname{\noexpand\ensuremath{\noexpand\mathscr{\@l}}}%
}

\newcommand{\spn}{\operatorname{span}}
\newcommand{\supp}{\operatorname{supp}}

\newcommand{\dist}{\operatorname{dist}}

\newcommand{\diam}{\operatorname{diam}}

\newcommand{\E}{\mathbb{E}}

\renewcommand{\d}{\operatorname{d}\!}

\newcommand{\bs}[1]{{\boldsymbol#1}}
\renewcommand{\d}{\operatorname{d\!}}

\newcommand{\isdef}{\mathrel{\mathrel{\mathop:}=}}

\definecolor{navy}{RGB}{102,153,255}
\definecolor{tuerkis}{RGB}{51,153,204}
\DeclareMathOperator*{\nnz}{\operatorname{nnz}}
\DeclareMathOperator*{\anz}{\operatorname{anz}}
\title[A fast direct solver for nonlocal operators]
{A fast direct solver for nonlocal operators in wavelet coordinates}
\author{Helmut Harbrecht}
\author{Michael Multerer}

\address{
Helmut Harbrecht,
Departement Mathematik und Informatik, Universit\"at Basel
Spiegelgasse~1, 4051 Basel, Schweiz
{\tt helmut.harbrecht@unibas.ch}
}

\address{
Michael Multerer,
Institute of Computational Science,
Universit\`a della Svizzera italiana,
Via Buffi~13, 6900 Lugano, Svizzera
{\tt michael.multerer@usi.ch}
}

\begin{document}
\begin{abstract}
In this article, we consider fast direct solvers for
nonlocal operators. The pivotal idea is to combine a
wavelet representation of the system matrix, yielding a
quasi-sparse matrix, with the nested dissection ordering scheme. 
The latter drastically reduces the fill-in during the factorization
of the system matrix by means of a Cholesky decomposition 
or an LU decomposition, respectively. This way, we end up with
the exact inverse of the compressed system
matrix with only a moderate increase of the number of nonzero entries 
in the matrix. 

To illustrate the efficacy of the approach, we conduct numerical
experiments for different highly relevant applications of nonlocal operators: 
We consider (i) the direct solution of boundary integral equations in three spatial 
dimensions, issuing from the polarizable continuum model, (ii) a parabolic 
problem for the fractional Laplacian in integral form and (iii) the fast simulation 
of Gaussian random fields.
\end{abstract}

\keywords{Nonlocal operator, Direct solver, 
Wavelet matrix compression,
Polarizable continuum model,
Fractional Laplacian,
Gaussian random fields}

\maketitle

\section{Introduction}
Various problems in science and engineering lead to
nonlocal operators and corresponding operator equations.
Examples arise from physical problems like field calculations
and Riesz energy problems, from machine learning, and also 
from stochastic simulations and uncertainty quantification.

Traditional discretizations of nonlocal operators result in densely 
populated system matrices. This feature renders the computation very costly in 
both respects, the computation time and computer memory requirements. 
Therefore, over recent decades, different ideas for the data sparse approximation 
of nonlocal operators have been developed. Most prominent examples of 
such methods are the fast multipole method \cite{GR}, the panel clustering 
\cite{HN}, the wavelet matrix compression \cite{BCR,DHS1}, and the hierarchical 
matrix format \cite{HA2}. These techniques are able to represent nonlocal 
operators in linear or almost linear cost with respect to the 
number of degrees of freedom used for their discretization.

The present article relies on a compression of the system matrix by
wavelets. Especially, the matrix representation of the nonlocal operator 
in wavelet coordinates is quasi-sparse, i.e.\ most matrix entries are 
negligible and can be treated as zero without compromising the overall
accuracy. Discarding the non-relevant matrix entries is called matrix 
compression. Roughly speaking, nonlocal operators become local operators 
in wavelet coordinates. A fully discrete version of the wavelet 
matrix compression has been developed in \cite{HS2}. It computes the 
compressed operator within discretization accuracy with linear cost.

Based on the sparsity pattern of the system matrix which is 
solely determined by the order of the underlying operator, we
employ a fill-in reducing reordering of the matrix entries by
means of nested dissection, see \cite{Geo73,LRT79}. This 
reordering in turn allows for the rapid inversion of the system 
matrix by the Cholesky decomposition or more generally by 
the LU decomposition. In particular, besides the rigorously
controllable error for the matrix compression in the wavelet 
format and the roundoff errors in the computation of the matrix
factorization, no additional approximation errors are introduced.
This is a major difference to other approaches for the discretization 
and the arithmetics of nonlocal operators, e.g.\ by means of hierarchical 
matrices. As the hierarchical matrix format is not closed under 
arithmetic operations, a recompression step after each arithmetic 
(block) operation has to be performed, which results into 
accumulating and hardly controllable consistency errors 
for matrix factorizations, see \cite{HA1,HA2}. 

In order to demonstrate the efficiency of the suggested approach, 
we consider applications from different fields. Namely, we consider 
(i) a boundary integral equation arising from the polarizable
continuum model in quantum chemistry as a classical example 
for a nonlocal operator equation, (ii) a parabolic problem for the 
fractional Laplacian, and finally (iii) the fast numerical simulation of 
Gaussian random fields as an important example from computational 
uncertainty quantification.

One of the most widespread methods to include solvent effects in 
quantum chemistry is by making use of a continuum dielectric: the 
solvent is represented by a continuum which surrounds the molecule.
Solute-solvent interactions are then described through appropriate 
functions supported on the molecule's surface. For an overview of 
continuum solvation models, we refer the reader to \cite{TMC}, and 
in particular for the polarizable continuum model to \cite{CM,CMT,MST}.
Wavelet matrix compression for the polarizable continuum model has been
considered in \cite{WaveletPCM2,WaveletPCM1}. Especially, in \cite{ILU}, 
the use of an incomplete Cholesky decomposition based preconditioner has 
been suggested, which is however inferior to the approach presented here. 

The fractional Laplacian is an operator which generalizes the notion of 
spatial derivatives to fractional orderss. It appears in image asnalysis, 
kinetic equations, phase transitions and nonlocal heat conduction, 
just to mention some applications. We refer to the review article
\cite{DGZ} and the references therein for further details. In particular,
we will focus here on the definition of the fractional Laplacian in its
integral form, as it can be found in \cite{EG} and also in \cite{DGZ}. 
To the best of our knowledge, the numerical treatment of the parabolic 
problem for the fractional Laplacian by means of wavelets has not 
been addressed in literature yet. 

The rapid simulation of (Gaussian) random fields with a prescribed 
covariance structure is of paramount importance in computational 
uncertainty quantification. The fast methods, which have been 
suggested so far are based on the computation of matrix square 
roots employing low-rank factorizations, block-wise low-rank 
factorizations, such as obtained by hierarchical matrices, or the 
discretization of the action of the matrix square root on a given 
vector by Krylov subspace methods. Other approaches compute the 
Karhunen-Lo\`eve expansion by circulant embeddings and fast Fourier 
techniques or employ the contour integral method. For more details on 
these methods, we refer to \cite{FKS18,GKNS+18,HPS15,HKS20,RH}.
In contrast to the previously mentioned approaches, we consider here 
the direct simulation of the random field by the Cholesky decomposition 
of the covariance matrix, which is very sparse in wavelet coordinates.

This article is organized as follows. Wavelet bases and 
their properties are specified in Section~\ref{sec:wavelets}.
Section~\ref{sec:WEM} briefly repeats the main features of the 
fully discrete wavelet matrix compression scheme from \cite{HS2}. 
Then, in Section~\ref{sect:ND}, for the sake of completeness,
the idea of nested dissection is briefly outlined.
Section~\ref{sct:applications} presents the three different 
applications considered in this article, while Section~\ref{sec:numerix} 
is devoted to related numerical experiments. Finally, Section~\ref{sec:conclusion} 
contains some concluding remarks.

In the following, in order to avoid the repeated use of generic but
unspecified constants, we write $C \lesssim D$ to indicate that $C$ 
can be bounded by a multiple of $D$, independently of parameters
which $C$ and $D$ may depend on. Then, $C \gtrsim D$
is defined as $D \lesssim C$, while we write $C \sim D$, iff
$C \lesssim D$ and $C \gtrsim D$.

\section{Wavelets and multiresolution Analysis}
	\label{sec:wavelets}
Let $D$ denote a domain in $\mathbb{R}^n$ or 
a manifold in $\mathbb{R}^{n+1}$. A 
\emph{multiresolution analysis} consists of a nested 
family of finite dimensional approximation spaces
\begin{equation}	\label{eq:hierarchy}
\{0\} = V_{-1} \subset V_0 \subset V_1 \subset \cdots 
\subset V_j\subset \cdots \subset L^2(D),
\end{equation}
such that
\[ \overline{\bigcup_{j\ge0} V_j} = L^2(D)\quad\text{and}\quad
\dim V_j\sim 2^{jn}.
\] We will
refer to \(j\) as the \emph{level} of \(V_j\) in the
multiresolution analysis.
Each space $V_j$ is endowed with a \emph{single-scale basis}
\[
\Phi_j = \{\varphi_{j,{\bs k}}:{\bs k}\in\Delta_j\},
\]
i.e.\ $V_j=\spn\Phi_j$, where $\Delta_j$ denotes a suitable 
index set with cardinality $|\Delta_j| \sim 2^{jn}$. For convenience, 
we shall in the sequel write bases on the form of row vectors, such 
that, for ${\bs v} = [v_k]_{k\in\Delta_j}\in\ell^2(\Delta_j)$, the 
corresponding function can simply be written as a dot product according
to
\[
v_j = \Phi_j{\bs v}=\sum_{k\in\Delta_j} v_k\varphi_{j,k}. 
\]
In addition, we shall assume that the single-scale 
bases $\Phi_j$ are \emph{uniformly stable}, this means that 
\[
\|{\bs v}\|_{\ell^2(\Delta_j)}
\sim \|\Phi_j{\bs v}\|_{L^2(D)}\quad\text{for all }
{\bs v}\in\ell^2(\Delta_j)
\]
uniformly in $j$, and that they satisfy the locality condition 
\[
\diam(\supp\varphi_{j,{\bs k}})\sim 2^{-j}.
\]

Additional properties of the spaces $V_j$ are required
for using them as trial spaces in a Galerkin scheme.
The approximation spaces shall have the \emph{regularity} 
\[
\gamma\isdef\sup\{s\in\mathbb{R}: V_j\subset H^s(D)\}
\]
and the \emph{approximation order} $d\in\mathbb{N}$, 
that is
\[
       d = \sup\Big\{s\in\mathbb{R}: \inf_{v_j\in V_j}\|v-v_j\|_{L^2(D)} 
       		\lesssim 2^{-js}\|v\|_{H^s(D)}\Big\}.
\]

Rather than using the multiresolution analysis corresponding to the 
hierarchy in \eqref{eq:hierarchy}, the pivotal idea 
of wavelets is to keep track of the increment of information 
between two consecutive levels $j-1$ and $j$. Since we have
$V_{j-1}\subset V_j$, we may decompose 
\[
V_j = V_{j-1}\oplus W_j,\ \text{i.e.}\ V_{j-1}\cup W_j = V_{j}\ \text{and}\  V_{j-1}\cap W_j = \{0\},
\]
with an appropriate \emph{detail space} $W_j$. Of practical interest 
is the particular choice of the basis of the detail space  $W_j$ in $V_j$. This basis 
will be denoted by
\[
  \Psi_j = \{\psi_{j,{\bs k}}:
  	{\bs k} \in \nabla_j\isdef \Delta_j\setminus \Delta_{j-1}\}.
\]
In particular, we shall assume that the collections $\Phi_{j-1}\cup\Psi_j$
form uniformly stable bases of $V_j$, as well. If $\Psi = 
\bigcup_{j\ge 0}\Psi_j$, where $\Psi_0\isdef \Phi_0$, is even a 
Riesz-basis of $L_2(D)$, then it is called a \emph{wavelet 
basis}. We require the functions $\psi_{j,{\bs k}}$ to be 
localized with respect to the corresponding level $j$, i.e.\ 
\[
\diam(\supp\psi_{j,{\bs k}}) \sim 2^{-j}, 
\]
and we normalize them such that 
\[
\|\psi_{j,{\bs k}}\|_{L_2(D)}\sim 1.
\]

At first glance it would be very convenient to deal with a single 
orthonormal system of wavelets. However, it has been shown
in \cite{DHS1,DPS4,S} that orthogonal wavelets are not optimal 
for the efficient approximation nonlocal operator equations. For this 
reason, we rather use \emph{biorthogonal wavelet bases}.
In this case, we also have a dual, multiresolution
analysis, i.e.\ dual single-scale bases and wavelets
\[
\widetilde{\Phi}_j = \{\widetilde{\varphi}_{j,{\bs k}}:{\bs k}\in\Delta_j\},\quad
\widetilde{\Psi}_j = \{\widetilde{\psi}_{j,{\bs k}}:{\bs k}\in\nabla_j\},
\] 
which are coupled to the primal ones by the orthogonality 
condition 
\[
(\Phi_j,\widetilde{\Phi}_j)_{L^2(D)} = {\bs I},\quad
(\Psi_j,\widetilde{\Psi}_j)_{L^2(D)} = {\bs I}. 
\]
The corresponding spaces $\widetilde{V}_j\isdef \spn\widetilde{\Phi}_j$
and $\widetilde{W}_j\isdef\spn\widetilde{\Psi}_j$ satisfy
\begin{equation}	\label{eq:space-coupling}
  V_{j-1}\perp \widetilde{W}_j, \quad \widetilde{V}_{j-1}\perp W_j.
\end{equation}
Moreover, the dual spaces are supposed to exhibit some approximation order
$\widetilde{d}\in\mathbb{N}$ and regularity $\widetilde{\gamma}>0$.

Denoting in complete analogy to the primal basis $\widetilde{\Psi} = 
\bigcup_{j\ge 0}\widetilde{\Psi}_j$, where $\widetilde{\Psi}_0 
\isdef\widetilde{\Phi}_0$, then every $v \in L^2(D)$ has
unique representations 
\[
v = \widetilde{\Psi} (v,\Psi)_{L^2(D)} = \Psi (v,\widetilde{\Psi})_{L^2(D)}
\]
such that
\[
 \|v\|_{L^2(D)}^2 \sim \sum_{j\ge 0}\sum_{k\in\nabla_j}
	\big\|(v,\widetilde\psi_{j,{\bs k}})_{L^2(D)}\big\|_{\ell^2(\nabla_j)}^2
 \sim \sum_{j\ge 0}\sum_{k\in\nabla_j}
	\big\|(v,\psi_{j,{\bs k}})_{L^2(D)}\big\|_{\ell^2(\nabla_j)}^2.
\]
In particular, relation \eqref{eq:space-coupling} implies that the 
wavelets exhibit {\em vanishing moments} of order $\widetilde{d}$, i.e.\
\begin{equation}	\label{eq:cancellation}
  \big|(v,\psi_{j,{\bs k}})_{L^2(D)}\big| \lesssim 2^{-j(1+\widetilde{d})}
	|v|_{W^{\widetilde{d},\infty}(\supp\psi_{j,{\bs k}})}.
\end{equation}
Herein, the quantity $|v|_{W^{\widetilde{d},\infty}(D)}
\isdef \sup_{|\boldsymbol\alpha|=
\widetilde{d}}\|\partial^{\boldsymbol\alpha}v\|_{L^\infty(D)}$
is the semi-norm in $W^{\widetilde{d},\infty}(D)$.
We refer to \cite{DA1} for further details.

Piecewise constant and bilinear wavelets which 
provide the above properties have been constructed in 
\cite{HS1,HS3}. In what follows, we will refer to the
wavelet basis of $V_J$ by $\Psi_J = \{\psi_\lambda:
\lambda\in\nabla_J\}$, where the multi-index 
$\lambda = (j,{\bs k})$ incorporates the scale 
$j=|\lambda|$ and the spatial location 
${\bs k} = {\bs k}(\lambda)$.

\section{Wavelet Matrix Compression}		\label{sec:WEM}
For a given domain or manifold $D$ and $q\in\mathbb{R}$, let
\[
\mathcal{A}:H^q(D)\to H^{-q}(D)
\]
denote a given (continuous and bijective) nonlocal operator of
order $2q$. According to the Schwartz kernel theorem, it can 
be represented in accordance with
\begin{equation}	 \label{eq:integral equation}
 (\mathcal{A}u)({\bs x}) = \int_D k({\bs x},{\bs y})
 	u({\bs y})\d{\bs y},
 \quad {\bs x}\in D,
\end{equation}
for a suitable kernel function $k\colon D\times D\to\mathbb{R}$.
The kernel functions under consideration are supposed 
to be smooth as functions in the variables ${\bs x}$ and 
${\bs y}$, apart from the diagonal $\{({\bs x},{\bs y})\in
D\times D: {\bs x}={\bs y}\}$ and may exhibit a singularity 
on the diagonal. Such kernel functions arise, for instance, from 
applying a boundary integral formulation to a second order elliptic 
problem \cite{SS,ST}. Typically, they decay like a negative power 
of the distance of the arguments which depends on the order 
$2q$ of the operator. More precisely, there holds
\begin{equation}	\label{eq:decay}
  \big|\partial_{\bs x}^{\boldsymbol\alpha}\partial_{\bs y}^{\boldsymbol\beta} 
  	k({\bs x},{\bs y})\big| \le c_{\boldsymbol\alpha,\boldsymbol\beta}
		\|{\bs x}-{\bs y}\|^{-n-2q-|\boldsymbol\alpha|-|\boldsymbol\beta|}.
\end{equation}
We emphasize that this estimate remains valid for the kernels of
arbitrary pseudodifferential operators, see \cite{DHS2} for the details.

Corresponding to the nonlocal operator from \eqref{eq:integral equation},
we may consider the operator equation \[\mathcal{A}u=f\] which gives rise
to the Galerkin approach:
\begin{align*}
&\text{find $u_J\in V_J$ such that}\\
&\quad (\mathcal{A}u_J,v_J)_{L^2(D)} = (f,v_J)_{L^2(D)}
\quad\text{for all $v_J \in V_J$}.
\end{align*}
Traditionally, this equation is discretized 
employing the single-scale basis of $V_J$ which results in densely
populated system matrices. If $N_J\sim 2^{Jn}$ denotes the 
number of basis functions in the space $V_J$, then the system 
matrix contains $\mathcal{O}(N_J^2)$ nonzero matrix entries. In 
contrast, by utilizing a wavelet basis in the Galerkin discretization,
we end up with a matrix that is quasi-sparse, i.e.\ it is compressible 
to $\mathcal{O}(N_J)$ nonzero matrix entries without compromising 
the overall accuracy. More precisely, by combining \eqref{eq:cancellation}
and \eqref{eq:decay}, we arrive at the decay estimate
\begin{equation}	\label{eq:decay estimate}
  (\mathcal{A}\psi_{\lambda'},\psi_\lambda)_{L^2(D)}
  	\lesssim \frac{2^{-(|\lambda|+|\lambda'|)(\widetilde{d}+n/2)}}
	{\operatorname{dist}(D_\lambda,D_{\lambda'})^{n+2q+2\widetilde{d}}}
\end{equation}
which is the foundation of the compression estimates in \cite{DHS1}. 
Herein, $D_\lambda\isdef\supp\psi_\lambda$ and 
$D_{\lambda'}\isdef\supp\psi_\lambda$ denote the 
convex hulls of the supports of the wavelets $\psi_\lambda$ and
$\psi_{\lambda'}$.

Based on \eqref{eq:decay estimate}, we shall neglect
all matrix entries for which the distance of 
the supports between the associated ansatz and test 
wavelets is larger than a level dependent cut-off 
parameter $\mathcal{B}_{j,j'}$. An additional compression, 
reflected by a cut-off parameter $\mathcal{B}_{j,j'}^s$, 
is achieved by neglecting several of those matrix entries, 
for which the corresponding trial and test functions 
have overlapping supports. 

To formulate this result, we introduce the abbreviation
$D_\lambda^s\isdef \operatorname{sing}\supp\psi_\lambda$
which denotes the {\em singular support} of the wavelet 
$\psi_\lambda$, i.e.\ that subset of $D$ where 
the wavelet is non-smooth.

\begin{theorem}[A-priori compression \cite{DHS1}]\label{thm:a-priori}
Let $D_\lambda$ and $D_\lambda^s$ be given as 
above and define the compressed system matrix ${\bs A}_J$, 
corresponding to the boundary integral operator $\mathcal{A}$, by
\begin{equation}		\label{eq:a-priori}
  [{\bs A}_J]_{\lambda,\lambda'}\isdef \begin{cases}
  \qquad \quad 0,& \dist(D_\lambda,D_{\lambda'}) > \mathcal{B}_{|\lambda|,|\lambda'|}\ \text{and $|\lambda|,|\lambda'|> 0$}, \\
  \qquad \quad 0,& \dist(D_\lambda,D_{\lambda'}) \le 2^{-\min\{|\lambda|,|\lambda'|\}}\ \text{and}\\
                 & \dist(D_\lambda^s,D_{\lambda'}) >
		 \mathcal{B}_{|\lambda|,|\lambda'|}^s\ \text{if $|\lambda'| > |\lambda|\geq 0$}, \\
		 & \dist(D_\lambda,D_{\lambda'}^s) >
		 \mathcal{B}_{|\lambda|,|\lambda'|}^s\ \text{if $|\lambda| > |\lambda'|\geq 0$}, \\
  (\mathcal{A}\psi_{\lambda'},\psi_\lambda)_{L^2(D)}, &\text{otherwise}. \end{cases}
\end{equation}
Fixing
\begin{equation}	\label{eq:parameters}
  a > 1, \qquad d < \delta < \widetilde{d} + 2q,
\end{equation}
the cut-off parameters $\mathcal{B}_{j,j'}$ and
$\mathcal{B}_{j,j'}^s$ are set according to
\begin{equation}	\label{eq:cut-off parameters}
  \begin{aligned}
  \mathcal{B}_{j,j'} \isdef a \phantom{'}\max\left\{ 2^{-\min\{j,j'\}},
	2^{\frac{2J(\delta-q)-(j+j')(\delta+\widetilde{d})}{2(\widetilde{d}+q)}}\right\}, \vspace*{2mm} \\
  \mathcal{B}_{j,j'}^s \isdef a \max\left\{ 2^{-\max\{j,j'\}},
	2^{\frac{2J(\delta-q)-(j+j')\delta-\max\{j,j'\}
	\widetilde{d}}{\widetilde{d}+2q}}\right\}.
  \end{aligned}
\end{equation}
Then, the system matrix ${\bs A}_J$ only has $\mathcal{O}(N_J)$
nonzero entries. In addition, the error estimate
\begin{equation}	\label{eq:errest}
    \|u-u_J\|_{H^{2q-d}(D)} \lesssim 2^{-2J(d-q)}\|u\|_{H^d(D)}
\end{equation}
holds for the solution $u_J$ of the compressed Galerkin system
provided that $u$ and $D$ are sufficiently regular.
\end{theorem}

\begin{figure}
\begin{center}
\begin{minipage}{0.49\textwidth}
\begin{center}
\includegraphics[scale=0.55,clip,trim=0 0 0 13,frame]{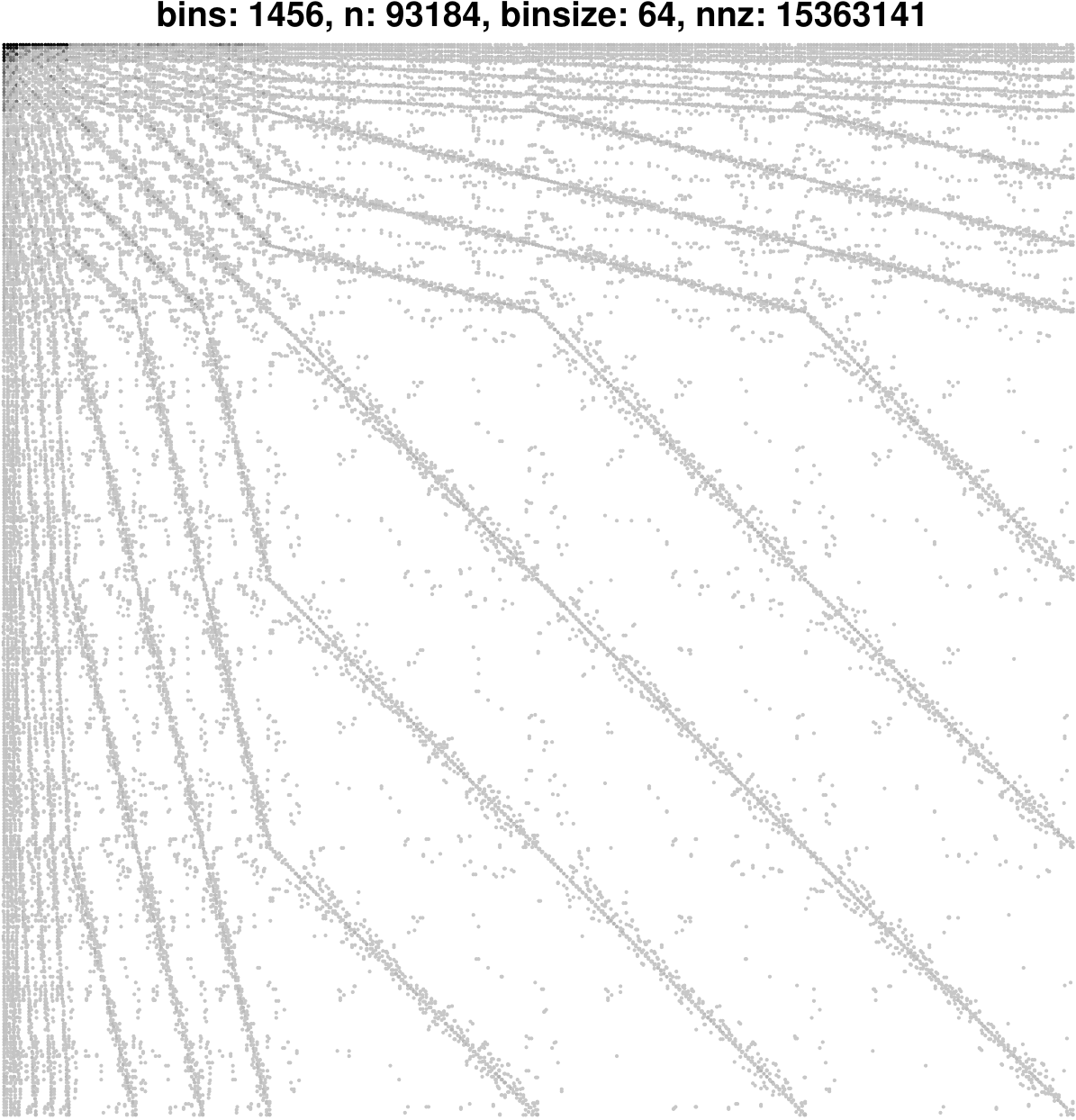}
$\nnz({\bs V}_J) = 15\,363\,141$
\end{center}
\end{minipage}
\begin{minipage}{0.49\textwidth}
\begin{center}
\includegraphics[scale=0.55,clip,trim=0 0 0 13,frame]{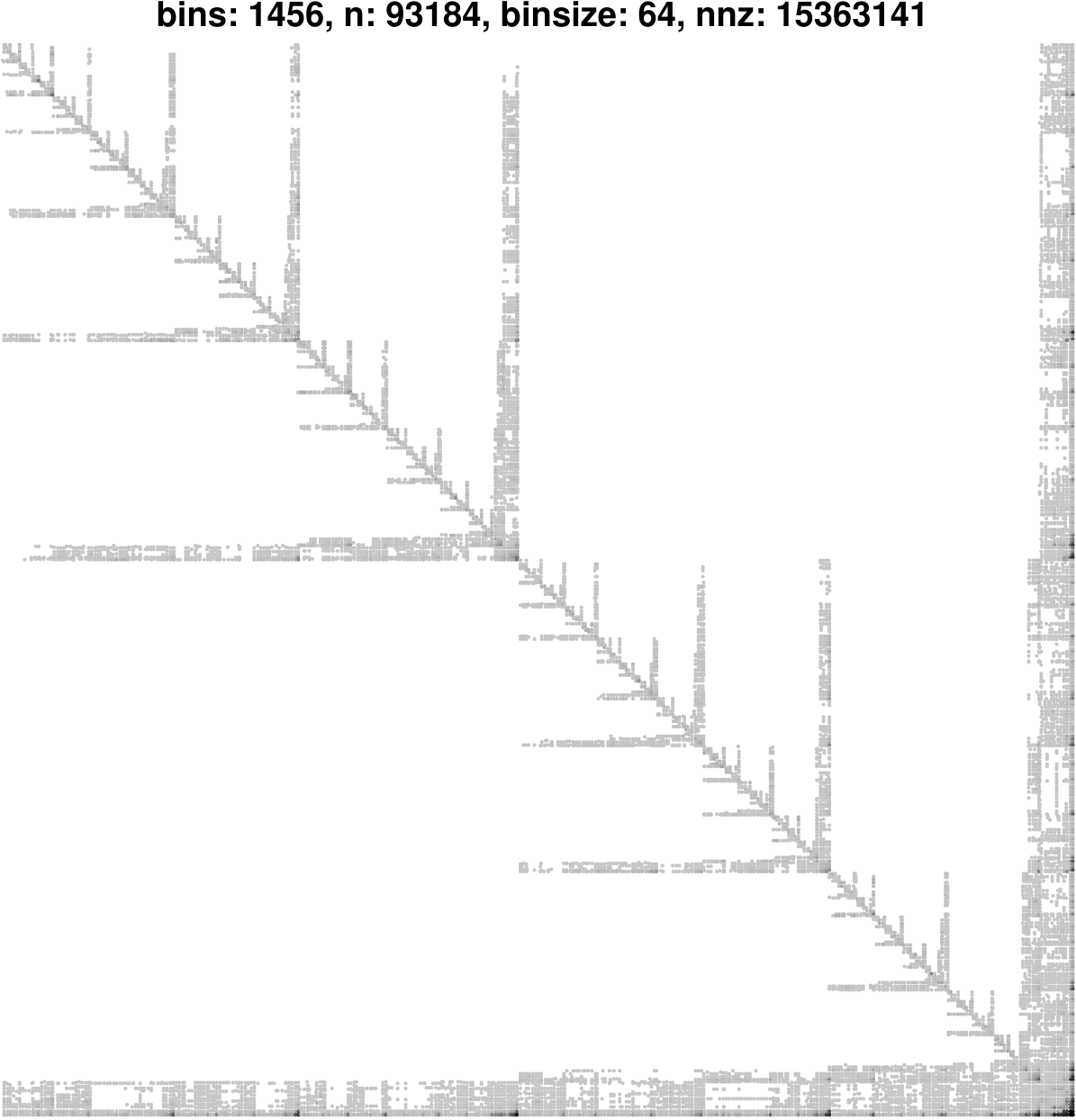}
$\nnz({\bs V}_{J,\text{ND}}) = 15\,363\,141$
\end{center}
\end{minipage}
\caption{\label{fig:compression}Sparsity patterns of \({\bs V}\) (left) and its nested
dissection reordering
\({\bs V}_{J,\text{ND}}\) (right) for the single layer operator on the benzene 
geometry and \(N_J=93184\). Each dot corresponds to a submatrix of size 
\(64\times 64\). Lighter blocks have less entries than darker blocks.}
\end{center}
\end{figure}

The compressed system matrix can be assembled with linear 
cost if the exponentially convergent 
$hp$--quadrature method proposed in \cite{HS2} is employed for
the computation of matrix entries. Moreover, 
for performing faster matrix-vector multiplications, an
additional a-posteriori compression might be applied
which reduces again the number of nonzero entries 
by a factor 2--5, see \cite{DHS1}. The pattern of the compressed 
system matrix shows the typical {\em finger structure}, see 
the left hand side of Figure~\ref{fig:compression}.

\section{Nested dissection}\label{sect:ND}
The representation of the system matrix corresponding to
a nonlocal operator with respect to an appropriate wavelet
basis leads to a quasi-sparse matrix, i.e.\ a matrix with many
small entries which can be neglected without compromising
accuracy. Performing a thresholding procedure as discussed 
in the previous section then yields a sparse system matrix
whose symmetric sparsity pattern is solely determined by the 
order of the underlying operator, see the left hand side of 
Figure~\ref{fig:compression}.

The factorization of the system matrix represented 
in the canoncial levelwise ordering leads to a massive 
fill-in. This means that a huge amount of nonzero 
entries is generated by a Cholesky decomposition or an LU decomposition, 
typically resulting in dense matrix factors. In order to obtain 
much sparser factorizations, we employ a nested dissection 
ordering, cf.\ \cite{Geo73,LRT79}, see the right hand side of 
Figure~\ref{fig:compression}.

Nested dissection is a divide and conquer algorithm whose foundation
is a graph theoretical observation. To each matrix
\({\bs A}\in\mathbb{R}^{N\times N}\) with a symmetric sparsity
pattern, we may assign an undirected graph \(G=(V,E)\) with
vertices \(V=\{1,2,\ldots,N\}\) and edges
\(E=\big\{\{i,j\}:a_{i,j}\neq 0\big\}\). Then, a symmetric permutation
\({\bs P}{\bs A}{\bs P}^\intercal\) of the rows and columns of \({\bs A}\)
amounts to a permutation \(\pi(V)\) of the nodes in \(V\).
In particular, we have the following important result from \cite{RTL76},
see also \cite{Par61}, which we formulate here only for the Cholesky decomposition
\({\bs P}{\bs A}{\bs P}^\intercal={\bs L}{\bs L}^\intercal\).

\begin{lemma}[\cite{RTL76}]
Assuming that no cancellation of nonzero entries
in the Cholesky decomposition
of \({\bs P}{\bs A}{\bs P}^\intercal\) takes place, then 
\(\ell_{i,j}\neq 0\) for \(i>j\), iff there is a path
\(i=v_1,v_2,\ldots, v_{k+1}=j\), \(k\geq 0\), in \(G\)
such that \(\pi(v_t)<\min\{\pi(i),\pi(j)\}\) for \(2\leq t\leq k\).
\end{lemma}

The lemma states that the Cholesky decomposition connects all
nodes \(i\) and \(j\), resulting in a nonzero entry \(\ell_{i,j}\),
for which there exists a path of nodes that have been eliminated
before \(i\) and \(j\).

Finding an optimal ordering is a hard problem in general. Therefore, we 
resort to the following strategy, which is known as nested dissection ordering:
We split \(V=V_1\cup V_2\cup S\)
such that \(E\cap\big\{\{v_1,v_2\}: v_1\in V_1,v_2\in V_2\big\}=\emptyset\),
i.e.\ the removal of the vertices of the \emph{separator} \(S\) and 
its adjacent edges results into two disjoint subgraphs.
Hence, employing an ordering which puts first the nodes into \(V_1\) and \(V_2\)
and afterwards the nodes in \(S\), leads to a matrix structure of the form
\[
{\bs P}{\bs A}{\bs P}^\intercal=
\begin{bmatrix}{\bs A}_{V_1,V_1} & & {\bs A}_{V_1,S}\\
 & {\bs A}_{V_2,V_2}& {\bs A}_{V_2,S}\\
{\bs A}_{S,V_1} & {\bs A}_{S,V_2}& {\bs A}_{S,S}\\\end{bmatrix}.
\]
Recursively applying this procudeure then yields a structure similar to
the one on the right hand side of Figure~\ref{fig:compression}.
For obvious reasons, it is desirable to have a minimal separator
\(S\), which evenly splits the \(V\) into two subsets, we refer to 
\cite{LRT79} and the references therein for a comprehensive 
discussion of this topic. In order to obtain suitable separators 
for the computations in this article, we will adopt the strategy 
from \cite{KK98}, which performs very well in terms of reducing 
the fill-in.

\section{Applications}\label{sct:applications}
\subsection{Polarizable continuum model}
Continuum solvation models are widely used to model
quantum effects of molecules in liquid solutions, compare
\cite{TMC} for an overview. In the {\em polarizable continuum 
model} (PCM), introduced in \cite{MST}, the molecule under 
study (the solute) is located inside a cavity $D$, surrounded 
by a homogeneous dielectric (the solvent) with dielectric 
constant $\epsilon\ge 1$. The solute-solvent interactions 
between the charge distributions which compose the solute 
and the dielectric are reduced to those of electrostatic origin.

For a given charge $\rho\in H^{-1}(D)$, located inside
the cavity, the solute-solvent interaction is expressed by the
{\em apparent surface charge} $\sigma\in H^{-1/2}(\partial D)$.
It is given by the integral equation
\begin{equation}	\label{eq:PCM}
 \mathcal{V}\sigma = \bigg(\frac{1+\epsilon}{2}+(1-\epsilon)\mathcal{K}\bigg)^{-1}
		\mathcal{N}_\rho-\mathcal{N}_\rho\quad\text{on $\partial D$},
\end{equation}
where $\mathcal{V}$ is the \emph{single layer operator}
\[
  (\mathcal{V}u)({\bs x}) = \int_{\partial D}
  	\frac{u({\bs y})}{4\pi\|{\bs x}-{\bs y}\|^3}\d\sigma_{\bs y},
\]
$\mathcal{K}$ is the \emph{double layer operator}
\[
  (\mathcal{K}u)({\bs x}) = \int_{\partial D} u({\bs y})
  	\frac{\langle{\bs n}({\bs y}),{\bs x}-{\bs y}\rangle}{4\pi\|{\bs x}-{\bs y}\|^3}\d\sigma_{\bs y},
\]
and $\mathcal{N}_\rho$ denotes the \emph{Newton potential} of the 
given charge
\[
  \mathcal{N}_\rho({\bs x})\isdef \int_{\partial D}\frac{\rho({\bs y})}{4\pi\|{\bs x}-{\bs y}\|}
  	\d{\bs y}.
\]

The discretization of the boundary integral equation
\eqref{eq:PCM} by means of a Galerkin scheme is as 
follows, compare \cite{HM1,HM2}: We make the ansatz 
\[
  \sigma_J = \sum_\lambda\sigma_\lambda\psi_\lambda
\]
and introduce the mass matrix
\[
{\bs G}_J = [(\psi_{\lambda'},\psi_\lambda)_{L^2(D)}]_{\lambda,\lambda'}
\]
and the system matrices
\[
 {\bs V}_J = [(\mathcal{V}\psi_{\lambda'},\psi_\lambda)_{L^2(D)}]_{\lambda,\lambda'},\quad
 {\bs K}_J = [(\mathcal{K}\psi_{\lambda'},\psi_\lambda)_{L^2(D)}]_{\lambda,\lambda'}.
\]
Then, for a given data vector ${\bs f}_J = [(\mathcal{N}_\rho,
\psi_\lambda)_{L^2(\partial D)}]_\lambda$, we need to solve the 
linear system of equations
\begin{equation}\label{eq:PCM2}
  {\bs V}_J\boldsymbol\sigma_J = {\bs G}_J
  	\bigg(\frac{1+\epsilon}{2}{\bs G}_J+(1-\epsilon) {\bs K}_J\bigg)^{-1}
  	{\bs f}_J-{\bs f}_J
\end{equation}
in order to determine the sought apparent surface charge.

In quantum chemical simulations, for example when 
solving the Hartree-Fock equations in a self consistent 
field approximation, one has to compute the interaction 
energies between the different particles. This amounts to 
the determination of different apparent surface charges. 
Therefore, the fast solution of \eqref{eq:PCM2} for multiple 
right hand sides is indispensable for fast simulations 
in quantum chemistry.

\subsection{Parabolic diffusion problem for the fractional Laplacian}
For a given domain $D\subset\mathbb{R}^n$ and $0<s<1/2$,
the fractional Laplacian $\mathcal{L}_s\colon H^s(D)/\mathbb{R}
\to H^{-s}(D)/\mathbb{R}$ is given by
\[
(\mathcal{L}_su)({\bs x})\isdef 2 \int_{D} \frac{u({\bs y}) - u({\bs x})}{\|{\bs x}-{\bs y}\|^{n+2s}}\d{\bs y},
	\quad {\bs x} \in D,
\]
compare \cite{EG,DGZ}. We intent to solve the following 
parabolic diffusion problem
\[
 \partial_t u - \mathcal{L}_s u = f\quad\text{in $D$}
\] 
for the fractional Laplacian.
To this end, we employ the $\theta$-scheme in time and a 
Galerkin discretization of the problem in space. This leads to the
linear system of equations
\begin{equation}\label{eq:theta}
 {\bs G}_J\frac{{\bs u}_J(t_{i+1})-{\bs u}_J(t_i)}{t_{i+1}-t_i} 
 - {\bs L}_J \big\{(1-\theta) {\bs u}_J(t_i) + \theta{\bs u}_J(t_{i+1})\big\}
 	= (1-\theta){\bs f}_J(t_i) + \theta{\bs f}_J(t_{i+1}).
\end{equation}
Here, 
\[
 {\bs L}_J = \big[(\mathcal{L}_s\psi_{\lambda'},
 	\psi_\lambda)_{L^2(D)}\big]_{\lambda,\lambda'},\quad
 {\bs G}_J = \big[(\psi_{\lambda'},\psi_\lambda)_{L^2(D)}\big]_{\lambda,\lambda'},
\]
are the system matrix of the fractional Laplacian and the
mass matrix, respectively. In each time step, we have hence to invert the
matrix ${\bs A}_J={\bs G}_J - \theta(t_{i+1}-t_i) {\bs L}_J$ for
computing the new solution ${\bs u}_J(t_{i+1})$ from the
solution  ${\bs u}_J(t_i)$ of the previous time step $t_i$.
A factorization of the system matrix \({\bs A}_J\)
is favourable in this situation since 
the system matrix does not change with time.

\subsection{Gaussian random fields}
Let \((\Omega,\mathcal{F},\mathbb{P})\) denote a complete
and separable probability space. 
We consider a Gaussian random field
\[
a\colon D\times\Omega\to\mathbb{R}
\]
with \emph{expectation} 
\[
\Ebb[a]({\bs x})\isdef\int_\Omega a({\bs x},\omega)\d\Pbb(\omega)
\] and \emph{covariance}
\[
\operatorname{Cov}[a]({\bs x},{\bs x}')\isdef\int_\Omega
\big(a({\bs x},\omega)-\Ebb[a]({\bs x})\big)
\big(a({\bs x}',\omega)-\Ebb[a]({\bs x}')\big)\d\Pbb(\omega).
\]
If the expectation and the covariance are known, we may represent
 \(a\) by its \emph{Karhunen-Lo\`eve expansion}
\[
a({\bs x},\omega)=\Ebb[a]({\bs x})+\sum_{k=1}^\infty\sqrt{\mu_k} a_k({\bs x})Y_k(\omega)
\]
Herein,  \((\mu_k,a_k)\), \(k=1,2,\ldots\),
denote the eigen pairs of the Hilbert-Schmidt operator
\[
(\mathcal{C}v)({\bs x})\isdef\int_D\operatorname{Cov}[a]({\bs x},{\bs x}')v({\bs x}')\d{\bs x}',
\]
while \(Y_1,Y_2,\ldots\) are independent and standard normally distributed random variables.

In order to discretize the Karhunen-Loeve expansion, we proceed in complete
analogy to \cite{HPS15} and compute the orthogonal projection of \(\mathcal{C}\) onto \(V_J\).
Let \({\bs{C}}_J\in\mathbb{R}^{N_J\times N_J}\) denote the corresponding coefficient matrix. Obviously, a suitable
basis for the orthogonal projection is given by the orthogonalized wavelet basis
\(\Psi_J{\bs G}^{-1/2}_J\), where
\[
{\bs G}_J = \big[(\psi_{\lambda'},\psi_\lambda)_{L^2(D)}\big]_{\lambda,\lambda'}
\]
denotes the mass matrix and \({\bs G}^{-1/2}_J\) is an inverse matrix root.
Thus, in accordance with \cite{HPS15}, we arrive at the discretized random field
\begin{equation}\label{eq:KLexp}
a_J({\bs x},\omega)=\Psi_J({\bs x}){\bs G}^{-1/2}_J\big(\overline{\bs a}_J + {\bs V}{\bs\Sigma}{\bs Y}(\omega)\big),
\end{equation}
where \({\bs Y}\isdef[Y_1,Y_2,\ldots, Y_{N_J}]^\intercal\) is a standard normally
distributed Gaussian vector, \(\overline{\bs a}_J\) is the orthogonal projection of
\(\E[a]\) onto \(V_J\) and \({\bs V}{\bs\Sigma}{\bs\Sigma}{\bs V}^\intercal={\bs C}_J\)
is the spectral decomposition of \({\bs C}_J\). 

Hence, for the Galerkin projection, we find
the identities
\[
\begin{aligned}
\overline{\bs a}^{G}_J&\isdef \big[(\E[a],\psi_{\lambda})_{L^2(D)}\big]_\lambda={\bs G}^{1/2}_J
\overline{\bs a}_J,\\
{\bs C}^G_J&\isdef\big[(\Ccal\psi_{\lambda'},\psi_{\lambda})_{L^2(D)}\big]_{\lambda,\lambda'}
={\bs G}^{1/2}_J{\bs C}_J{\bs G}^{1/2}_J.
\end{aligned}\] 
Hence, we obtain for the expectation
\[
\overline{\bs a}_J={\bs G}^{-1/2}_J\overline{\bs a}^{G}_J,
\]
while the covariance satisfies
\[
{\bs C}_J =({\bs V}{\bs\Sigma})({\bs V}{\bs\Sigma})^\intercal= {\bs G}^{-1/2}_J{\bs C}^G_J{\bs G}_J^{-1/2}=({\bs G}^{-1/2}_J{\bs L})({\bs G}^{-1/2}_J{\bs L})^\intercal,
\]
where \({\bs L}{\bs L}^\intercal={\bs C}^G_J\) is the Cholesky decomposition of \({\bs C}^G_J\).

It is
well known that any two matrix square roots
\({\bs R}{\bs R}^\intercal=\widehat{\bs R}\widehat{\bs R}^\intercal={\bs C}_J\)
only differ by an isometry. Consequently,
there exists an orthogonal matrix \({\bs Q}\)
such that \({\bs V}{\bs\Sigma}={\bs G}^{-1/2}_J{\bs L}{\bs Q}\).
Introducing the transformed random vector 
\(
{\bs X}(\omega)\isdef{\bs Q}{\bs Y}(\omega)
\), we infer
\[
\int_\Omega{\bs X}(\omega){\bs X}^\intercal({\omega})\d\mathbb{P}(\omega)
={\bs Q}\int_\Omega{\bs X}(\omega){\bs X}^\intercal({\omega})\d\mathbb{P}(\omega){\bs Q}^\intercal={\bs Q}{\bs I}{\bs Q}^\intercal={\bs I}.
\]
Therefore, \({\bs X}(\omega)\) is a standard normally
distributed Gaussian vector as well and we end up with the representation
\begin{equation}\label{eq:CholFastField}
a_J({\bs x},\omega)=\Psi_J({\bs x}){\bs G}^{-1}_J
\big(\overline{\bs a}^{G}_J + {\bs L}{\bs X}(\omega)\big).
\end{equation}

If the singular values in the Karhunen-Loeve expansion \eqref{eq:KLexp} 
decay only slowly, as it is typically the case for non-smooth
covariance functions, the numerical solution of the associated
eigenvalue problem becomes prohibitive. In such cases, the direct 
simulation of the random field by means of a sparse Cholesky 
decomposition in \eqref{eq:CholFastField} is computationally 
superior. Therefore, the proposed method is in particular useful for
rough random fields, which issue from non-smooth covariance 
functions. 

We remark that in case of a smooth covariance function, the 
matrix \({\bs C}_J^G\) becomes numerically positive semidefinite.
In such cases, the pivoted Cholesky decomposition, see \cite{HPS}, 
immediately yields a low-rank decomposition of \({\bs C}_J^G\) and 
is computationally superior to other data sparse representations,
we refer to \cite{HPS15} for a comprehensive discussion of this 
complementary approach.

\section{Numerical Results}\label{sec:numerix}
In order to obtain consistent computation times,
all computations reported in this section have been
carried out on a single core of a compute server with 
Intel Xeon E5-2650 v3 @2.30GHz CPUs and 512GB DDR4 
@2133MHz main memory. The implementation of the wavelet matrix
compression has been performed in ANSI C, while nested dissection
orderings, Cholesky decompositions and LU decompositions
have been computed using \textsc{Matlab} 2020a, compare 
\cite{MLB20}.

\subsection{Polarizable continuum model}
\begin{figure}[htb]
\begin{center}
\includegraphics[scale=0.16,clip,trim=290 30 280 110]{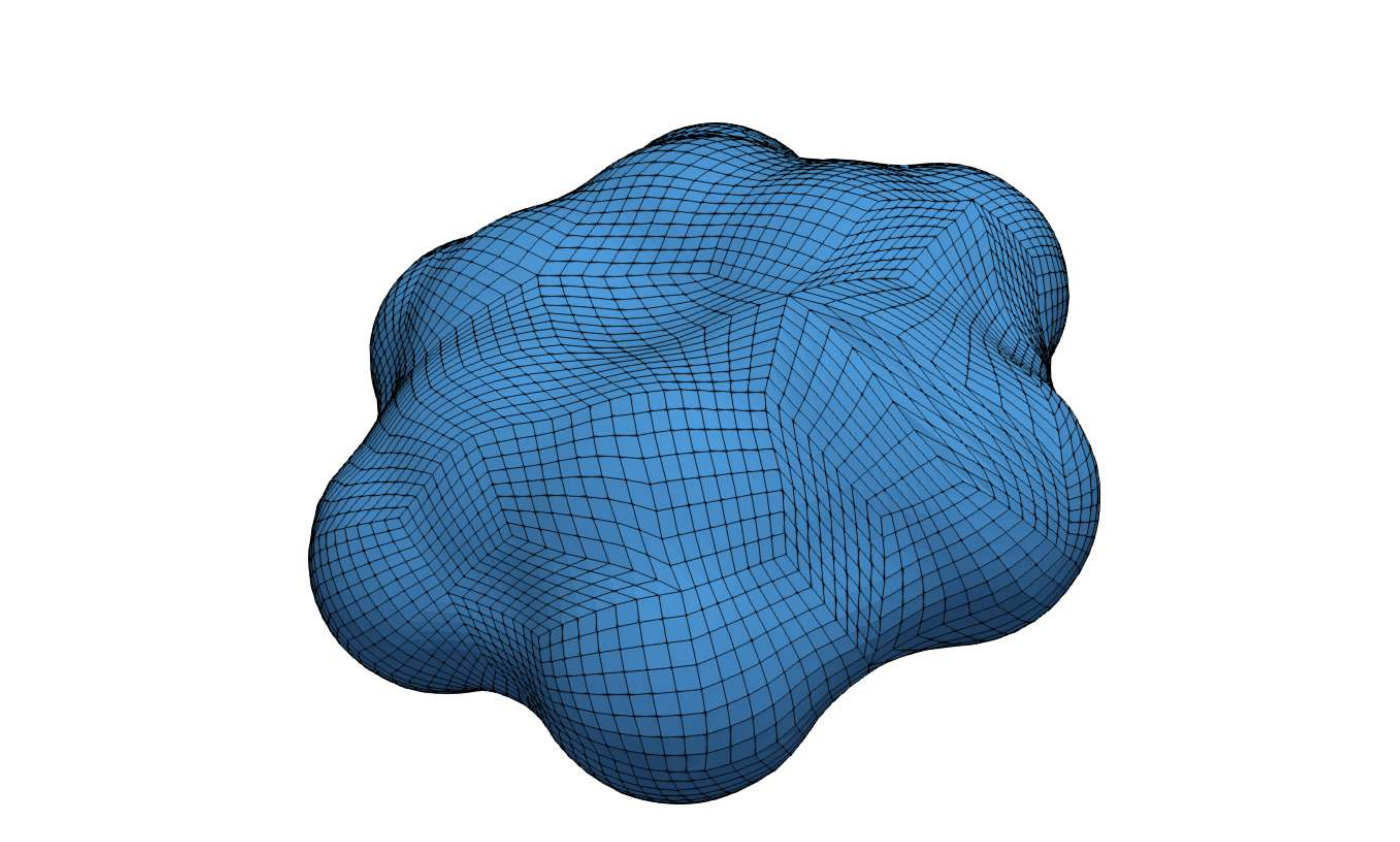}
\caption{\label{fig:benzene}Solvent excluding surface of benzene.}
\end{center}
\end{figure}

For PCM, we consider the assembly and the factorization of the
matrices
\[
  {\bs V}_J\quad\text{and}\quad{\bs A}_J\isdef
  	\frac{1+\epsilon}{2}{\bs G}_J+(1-\epsilon) {\bs K}_J,
\]
which are required to set up the linear system
of equations \eqref{eq:PCM2}. The dielectric constant is chosen as \(\epsilon=78.39\), 
corresponding to the solvent water. The molecule under study is 
benzene, whose solvent excluding surface is depicted in Figure
\ref{fig:benzene}. For the wavelet matrix compression, we use 
piecewise constant wavelets with three vanishing moments, 
as developed in \cite{HS1}. 

\begin{table}[htb]
\begin{center}
\begin{tabular}{|r|r|r|r|r|r|r|}\hline
$N_J$ & $t_\text{WEM}$ &  $t_\text{ND}$ &   $t_\text{Chol}({\bs V}_J)$ &
$t_\text{LU}({\bs A}_J)$ & $\anz({\bs V}_J)$  & $\anz({\bs L}_J)$ \\\hline
1456  &   0.96  &  0.032    &  0.19  &  0.22  & 134  &  209\\
5824  &   7.10  &  0.15    &   0.19  &  1.59  & 163  &  328\\
23296  &  52.32  &  0.74  &  1.05  & 14.03  & 159  &  419 \\
93184  &  309.35 &   3.91   &  7.66  & 189.70  & 165  &  639\\
372736  &  2265.91  &  19.38   & 58.28  & 2211.57  & 173  &  890\\
\hline
\end{tabular}
\caption{\label{tab:PCM}Computation times and 
numbers of nonzero entries in case of PCM.}
\end{center}
\end{table}

Table~\ref{tab:PCM} shows the numerical results. The first 
column labelled \(N_J\) corresponds to the number of surface
elements. The second column labelled \(t_\text{WEM}\) contains
the combined computation times in seconds for the assembly of \({\bs V}_J\)
and \({\bs A}_J\) by the wavelet matrix compression. The third column, 
labelled by \(t_\text{ND}\), provides the computation times in 
seconds for the nested dissection ordering, as \({\bs V}_J\)
and \({\bs A}_J\) have identical sparsity patterns, the same
reordering can be applied to both of them. The fourth
column labelled by \(t_\text{Chol}({\bs V}_J)\) denotes
the computation times in seconds for the Cholesky factorization 
of \({\bs V}_J\), while the fifth column shows the computation 
times in seconds of the LU decomposition of \({\bs A}_J\). For 
the purpose of measuring the fill-in relative to the matrix size, we
introduce for a sparse matrix \({\bs A}\in\mathbb{R}^{N\times N}\)
the \emph{average number of nonzeros per row}
\[
\anz({\bs A})\isdef \frac{\nnz({\bs A})}{N}.
\]
The values for \(\anz({\bs V}_J)\) and \(\anz({\bs A}_J)\) are 
identical and can be found in the sixth row of Table~\ref{tab:PCM},
while the values \(\anz({\bs L}_J)\) for the Cholesky factor of 
\(\anz({\bs V}_J)\) are given in the last column of the table. The 
sparsity pattern of the
L-factor of the LU decomposition of \({\bs A}_J\) coincides with that
of ${\bs L}_J$, while the U-factor has somewhat less coefficients compared 
to the L-factor. Hence, the average number of nonzeros per row for 
the LU factorization of \({\bs A}_J\) matches that of \(\anz({\bs L}_J)\). 

\begin{figure}[htb]
\begin{center}
\setlength{\fboxrule}{0.2pt}
\begin{minipage}{0.31\textwidth}
\begin{center}
\includegraphics[scale=0.4,clip,trim=0 0 0 13,frame]{Spattern_PCM}
$\nnz({\bs V}_J) = 15\,363\,141$
\end{center}
\end{minipage}
\begin{minipage}{0.31\textwidth}
\begin{center}
\includegraphics[scale=0.4,clip,trim=0 0 0 13,frame]{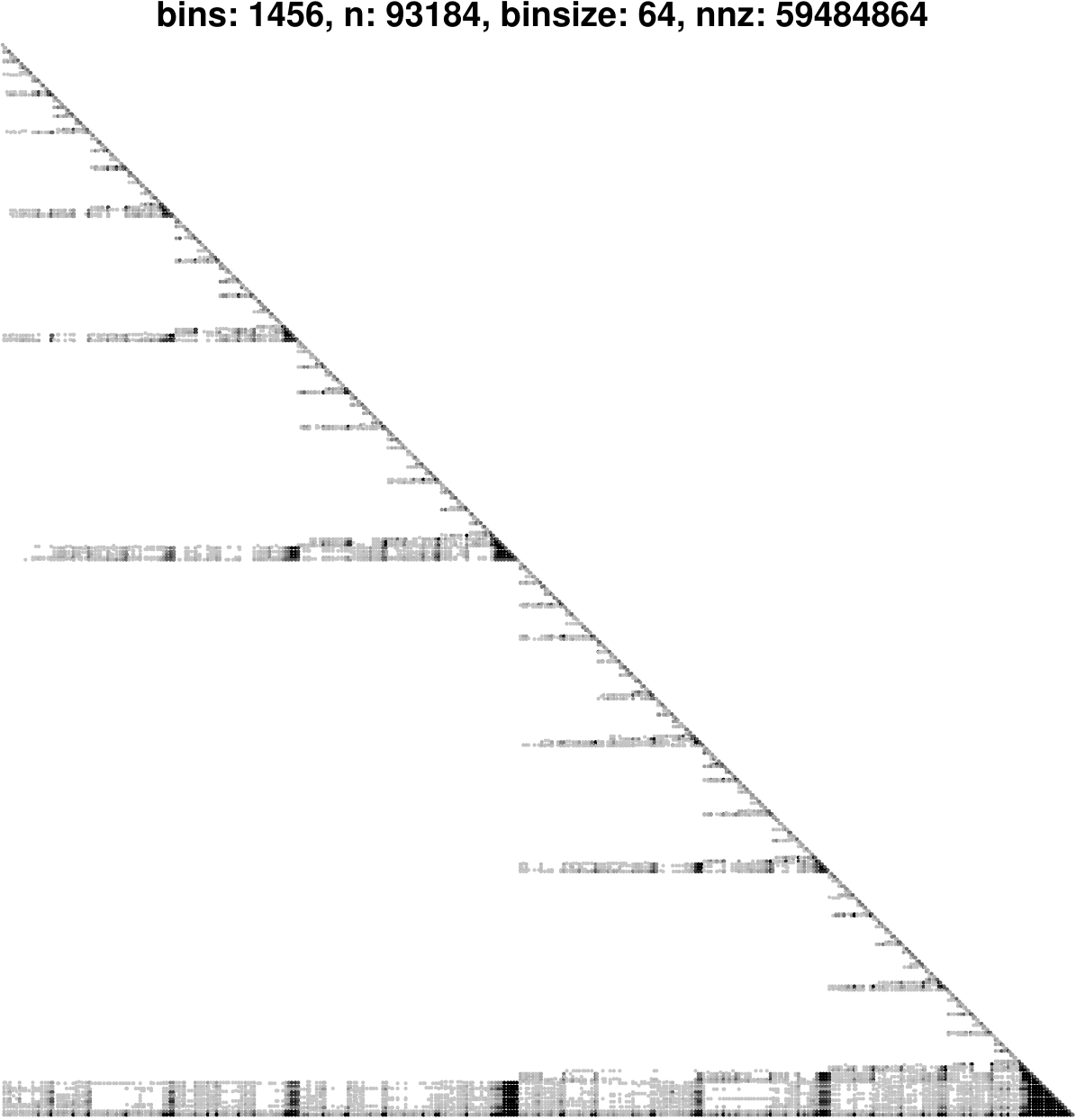}
$\nnz({\bs L}_J) = 59\,484\,864$
\end{center}
\end{minipage}
\begin{minipage}{0.31\textwidth}
\begin{center}
\includegraphics[scale=0.4,clip,trim=0 0 0 13,frame]{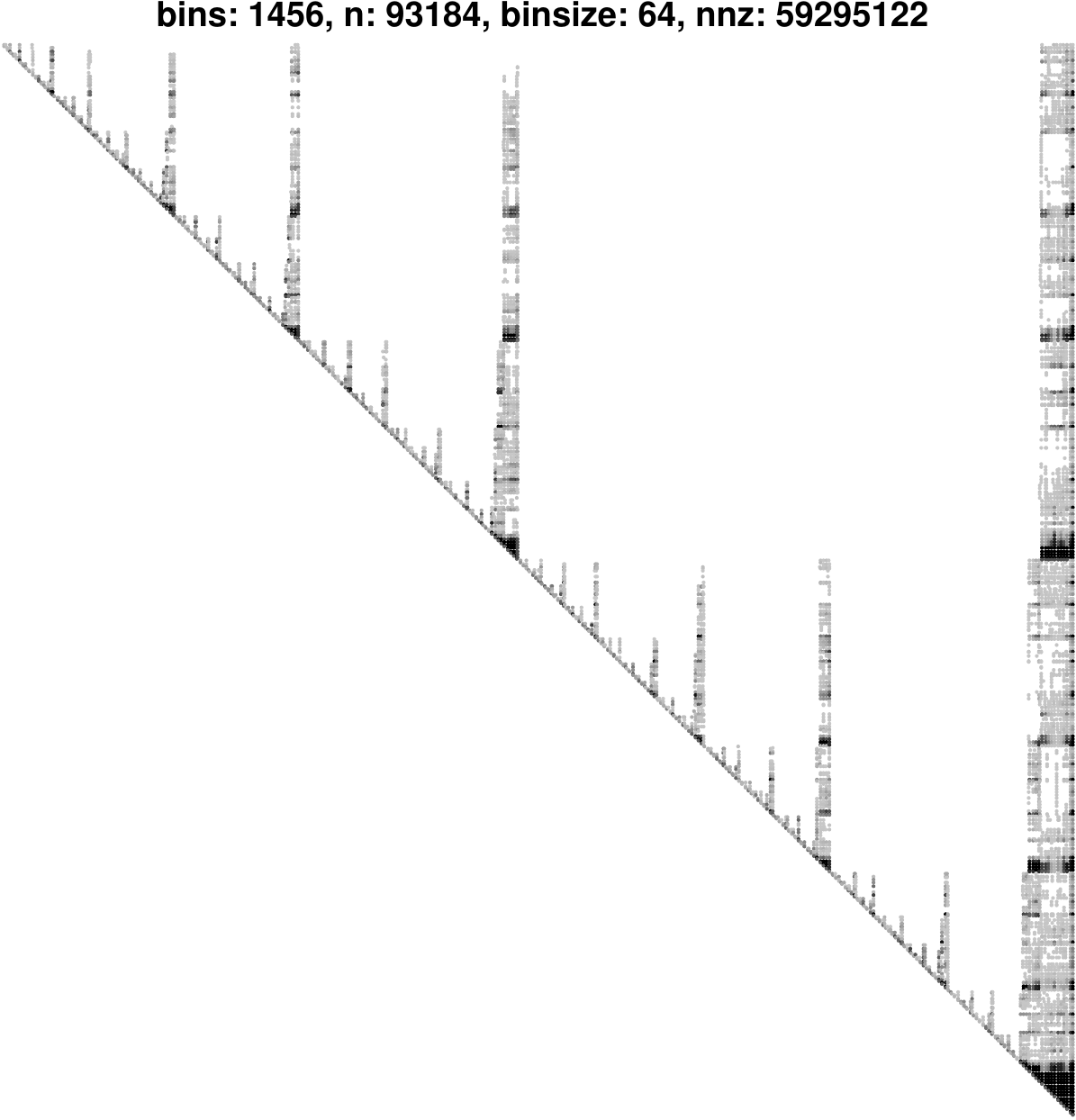}
$\nnz({\bs U}_J) = 59\,295\,122$
\end{center}
\end{minipage}
\caption{\label{fig:PCM}Sparsity patterns of \({\bs V}_J\) as well as \({\bs A}_J\) (left), 
the Cholesky factor \({\bs L}_J\) (middle), and the upper triangular matrix \({\bs U}_J\)
(right) for PCM and \(N_J=93\,184\). Each dot corresponds to a submatrix of size 
\(64\times 64\). Lighter blocks have less entries than darker blocks. }
\end{center}
\end{figure}

The sparsity patterns for \({\bs V}_J\), \({\bs L}_J\) and for the U-factor 
are shown in Figure~\ref{fig:PCM} for \(N_J=93\,184\). In order 
to obtain a neat representation in this figure, we have merged
matrix blocks of size \(64\times 64\). Darker blocks have a higher
density of entries, while lower blocks have a lower density of entries. 

As can clearly be inferred from Table~\ref{tab:PCM}, the 
times for the computing the nested dissection ordering and the 
subsequent Cholesky decomposition are negligible with respect 
to the wavelet matrix compression. The average number of nonzero 
entries per row stays rather low and only increases for increasing 
numbers of unknowns up to a factor of approximately 5 for \(N_J=372\,736\).

It can be seen from the fifth column of Table~\ref{tab:PCM} that
the LU decomposition is significantly slower than the Cholesky
decomposition. This issues from the fact that all matrices are 
stored in a sparse column major format, resulting in a large 
overhead for the access of matrix rows. We remark that, in 
principle, the LU decomposition could be accelerated by an 
appropriate data structure, which enables direct access 
to the rows of the matrix as well.
 
\subsection{Parabolic diffusion problem for the fractional Laplacian}
We consider the parabolic problem for the fractional Laplacian with 
\(s=\nicefrac{3}{8}\). The right hand side is a Gaussian heat
spot moving on a circular trajectory, given by
\[
f({\bs x},t)=100\exp\Big(-40\big(x_1-\cos(2\pi t)\big)^2-40\big(x_2-\sin(2\pi t)\big)^2\Big),
\]
while the initial condition is set to 0.
For the solution of the ordinary differential equation in time, 
we employ the \(\theta\)-scheme with \(\theta=\nicefrac{1}{2}\), which 
yields the Crank-Nicolson method \cite{CN}. The time interval 
is given by \([0,T]=[0,3]\) and we discretize this time interval
by 150 equidistant time steps. 

\begin{figure}[htb]
\begin{center}
\includegraphics[scale=0.14,clip,trim=510 160 510 160]{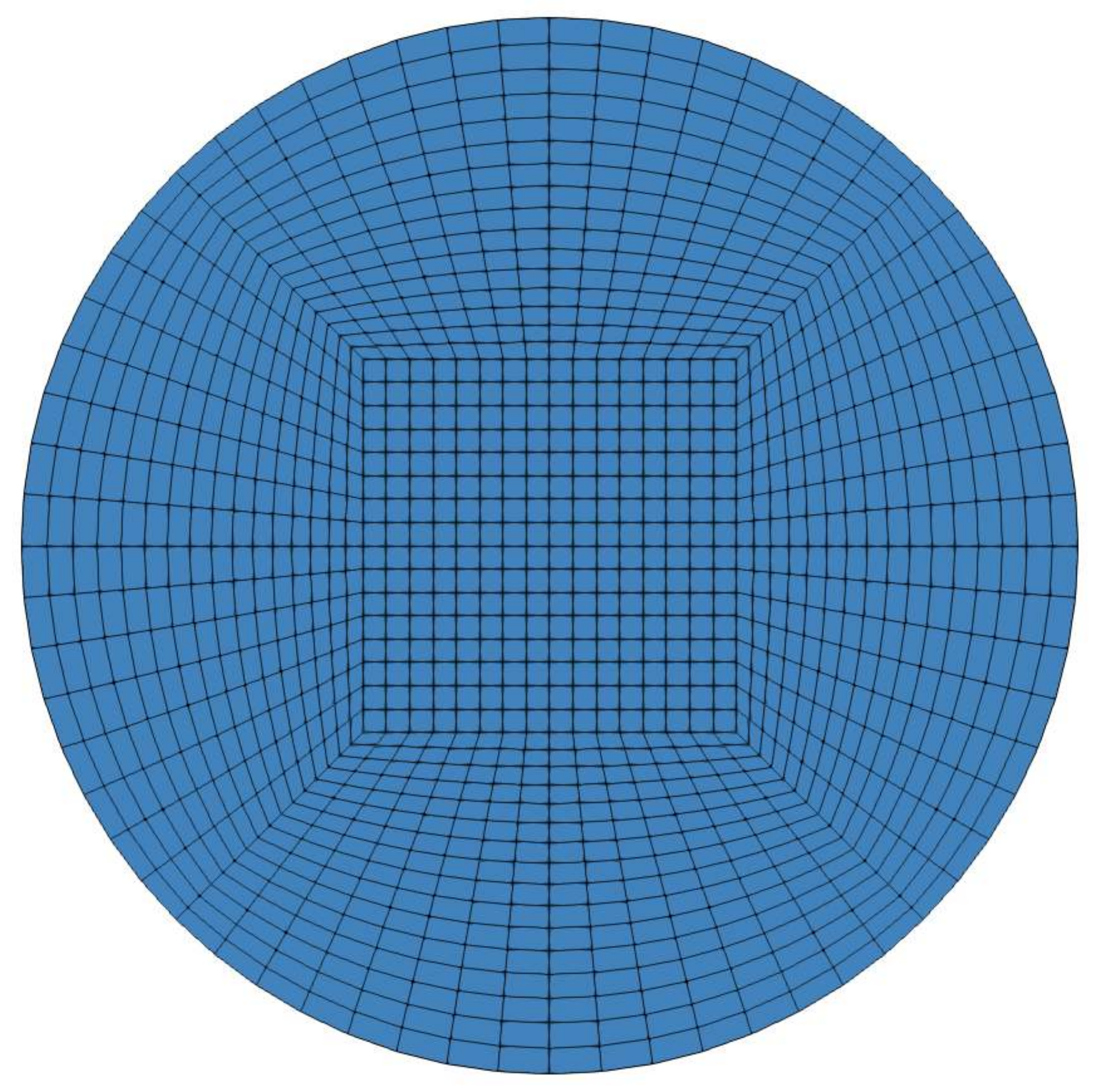}
\caption{\label{fig:disc}Computational geometry for the fractional Laplacian.}
\end{center}
\end{figure}

\begin{figure}[htb]
\begin{center}
\begin{tikzpicture}
\draw(0,0) node{\includegraphics[scale=0.1,clip,trim=1050 120 100 120]{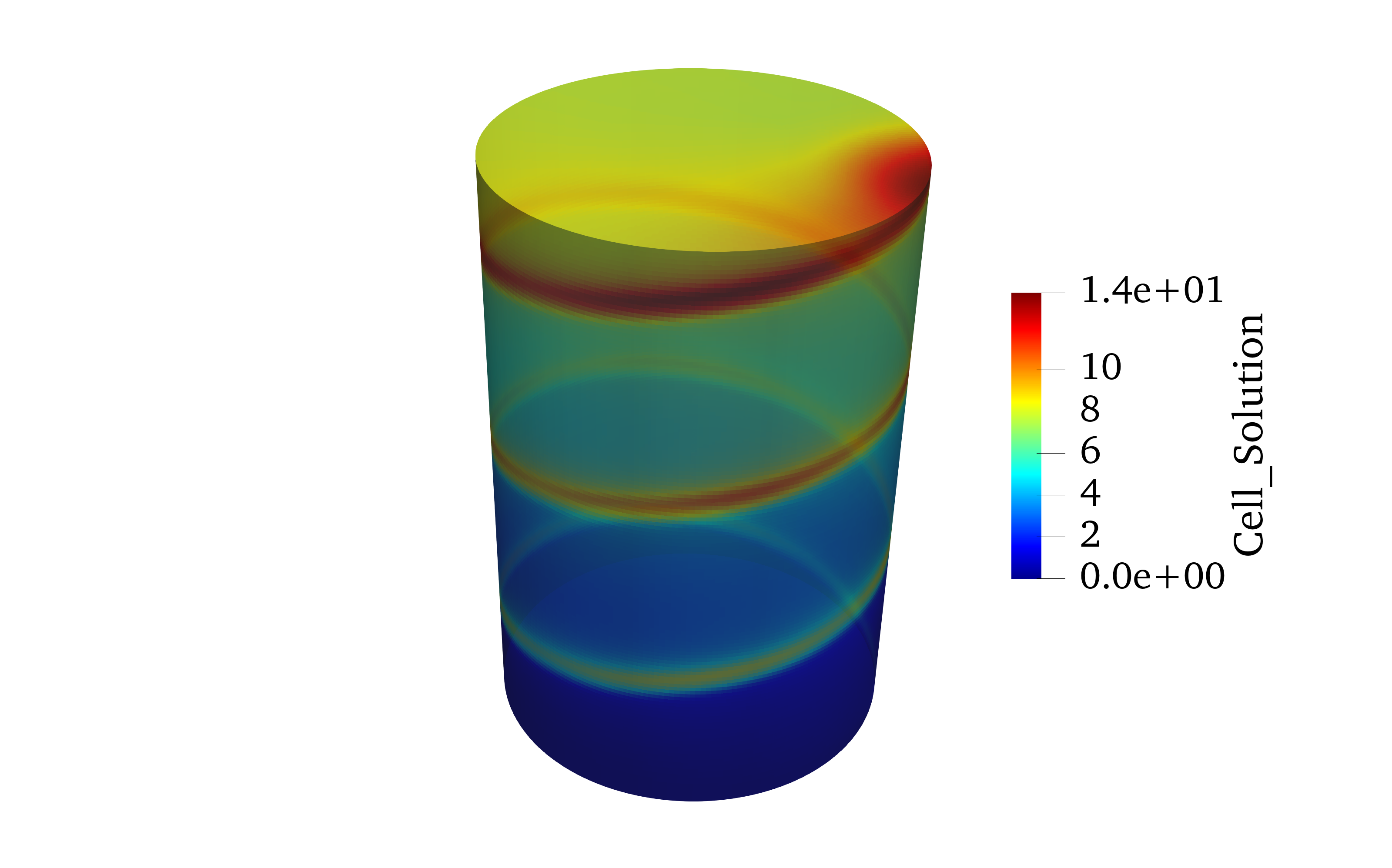}};
\fill[color=white] (1.3,-2) rectangle (3,2);
\draw(1.45,1.2) node {$14$};
\draw(1.45,0.6) node {$10$};
\draw(1.45,0.25) node {$\phantom{0}8$};
\draw(1.45,-0.1) node {$\phantom{0}6$};
\draw(1.45,-0.45) node {$\phantom{0}4$};
\draw(1.45,-0.8) node {$\phantom{0}2$};
\draw(1.45,-1.15) node {$\phantom{0}0$};
\draw(2.4,0) node[rotate=90] {Temperature};

\end{tikzpicture}
\caption{\label{fig:solFL}Solution of the fractional Laplacian on the space-time
cylinder.}
\end{center}
\end{figure}
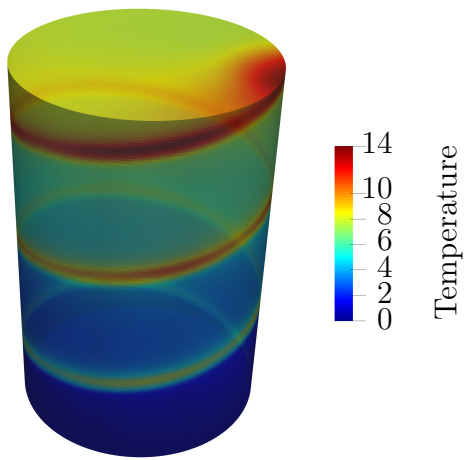

In the \(\theta\)-scheme \eqref{eq:theta}, we need to invert the matrix
\[
{\bs A}_J={\bs G}_J + \theta(t_{i+1}-t_i) {\bs L}_J,
\]
which is assembled with the help of wavelet matrix compression
using Haar wavelets. The computational geometry is the unit disc 
depicted in Figure~\ref{fig:disc}. The solution on the space-time 
cylinder is visualized for \(N_J= 20\,480\) in Figure~\ref{fig:solFL}.

\begin{table}[htb]
\begin{center}
\begin{tabular}{|r|r|r|r|r|r|r|}\hline
$N_J$ & $t_\text{WEM}$ &  $t_\text{ND}$ & $t_\text{Chol}({\bs A}_J)$ & 
\(t_{\theta}\)  & $\anz({\bs A}_J)$  & $\anz({\bs L}_J)$ \\\hline
   1280 & 0.63       &  0.02   &   0.01  & 0.19     &      88   &  140\\
   5120 &  3.35      &  0.11   &   0.10  & 0.96     &     106   &  222\\
  20480 &  22.20     &  0.53   &   0.59  & 5.34     &     117   &  315\\
  81920 &  161.83    &  2.71   &   3.89  & 27.94    &     129   &  474\\
 327680 &   724.60   & 12.87   &  23.49  & 131.30   &     135   &  576\\
1310720 &  3875.30   & 62.62   & 161.65  & 633.50   &     142   &  743\\
\hline
\end{tabular}
\caption{\label{tab:FLap}Computation times and numbers of 
nonzero entries in case of the fractional Laplacian.}
\end{center}
\end{table}

\begin{figure}[htb]
\begin{center}
\setlength{\fboxrule}{0.2pt}
\begin{minipage}{0.49\textwidth}
\begin{center}
\includegraphics[scale=0.55,clip,trim=0 0 0 13,frame]{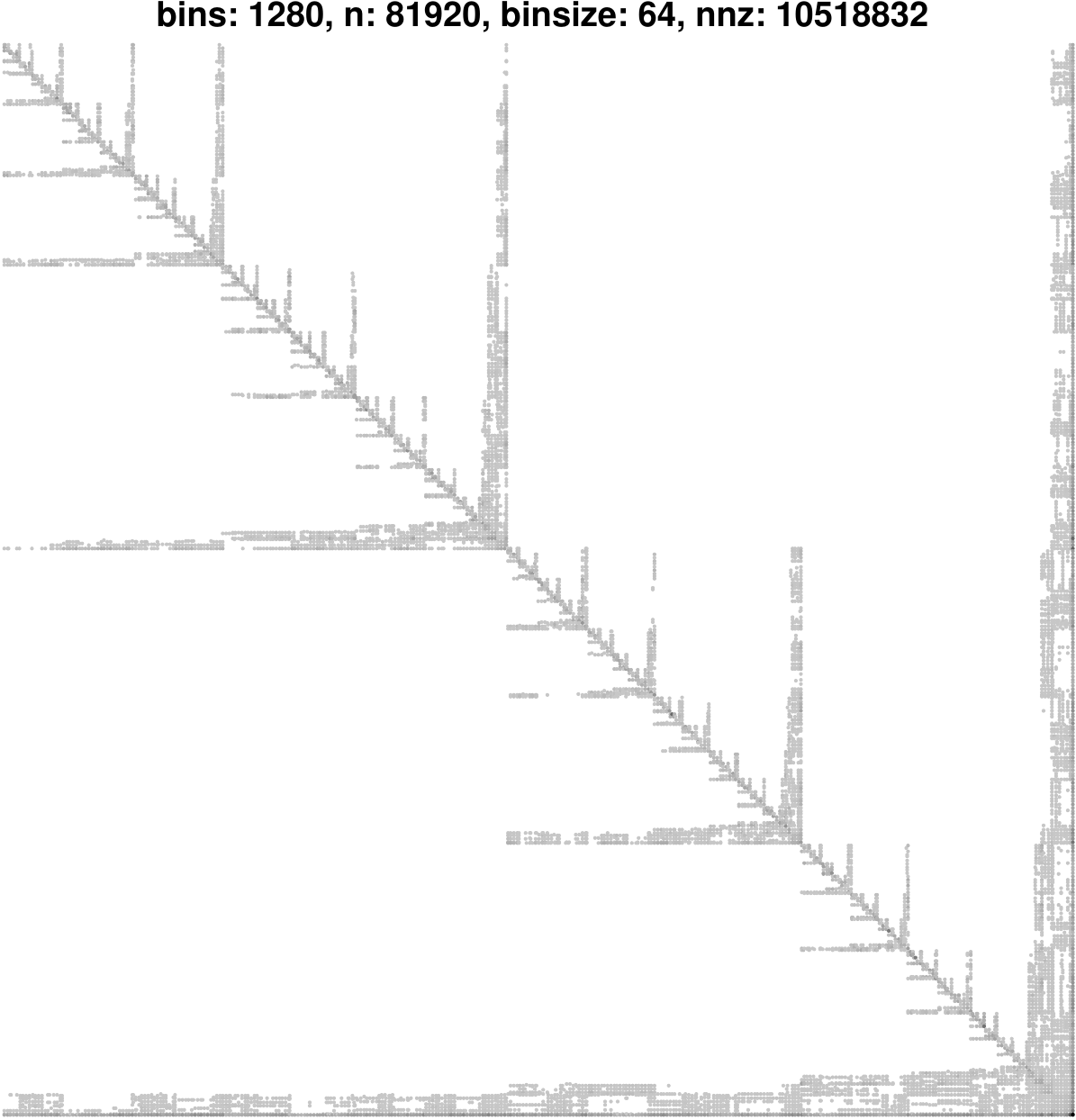}
$\nnz({\bs A}_J) = 10\,518\,832$
\end{center}
\end{minipage}
\begin{minipage}{0.49\textwidth}
\begin{center}
\includegraphics[scale=0.55,clip,trim=0 0 0 13,frame]{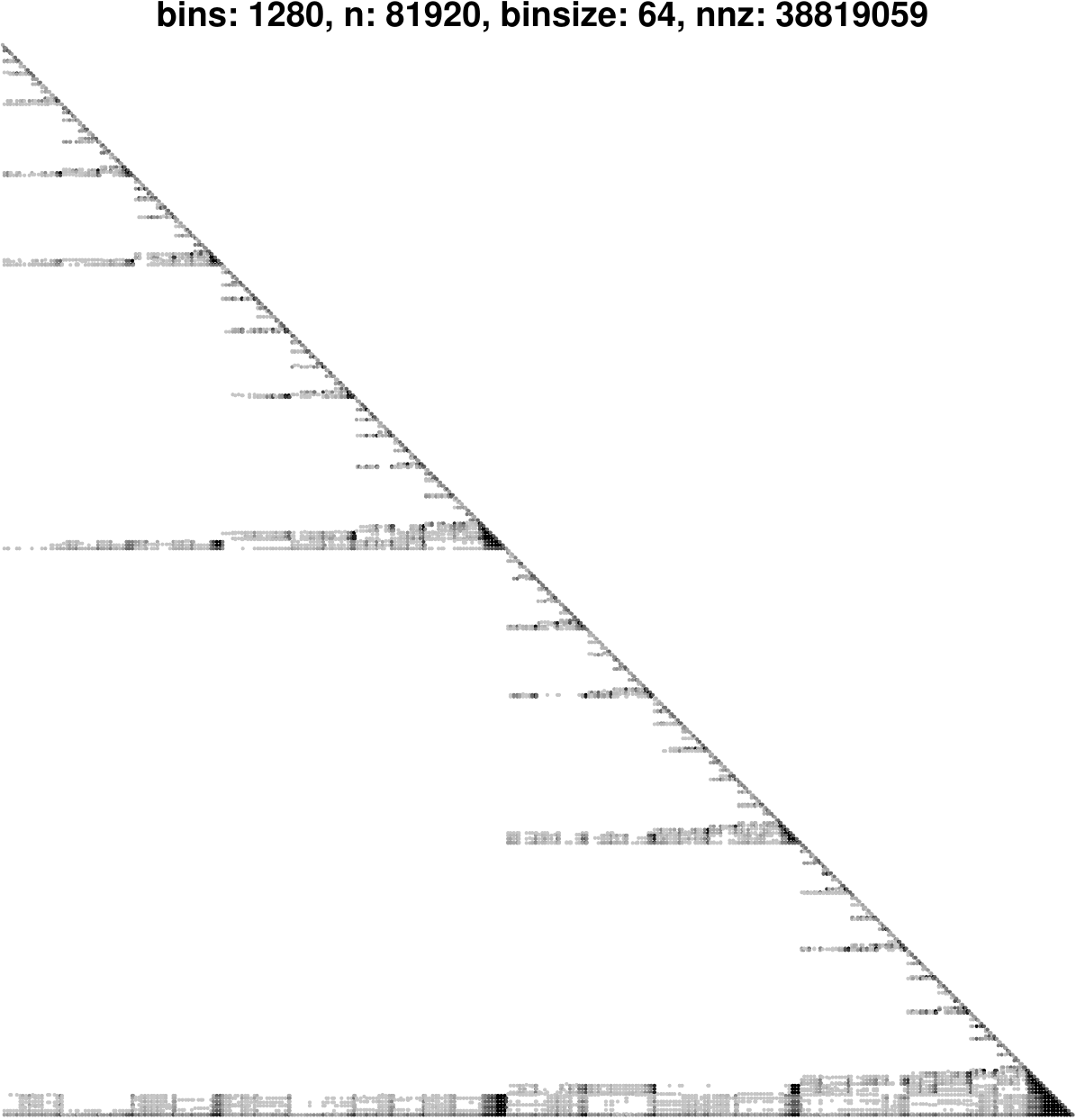}
$\nnz({\bs L}_J) = 38\,819\,059$
\end{center}
\end{minipage}
\caption{\label{fig:FractionalS}Sparsity patterns of \({\bs A}_J\) (left) and the Cholesky factor
\({\bs L}_J\) (right) for the fractional Laplacian and \(N_J=81\,920\). Each dot corresponds to 
a submatrix of size \(64\times 64\). Lighter blocks have less entries than darker blocks.}
\end{center}
\end{figure}

We have tabulated the similar characteristics from the previous 
example for the matrix \({\bs A}_J\) in Table~\ref{tab:FLap}, while
the sparsity patterns of \({\bs A}_J\) and \({\bs L}_J\) for \(N_J=81\,920\)
after reordering are shown in Figure~\ref{fig:FractionalS}. We observe 
a similar behaviour as for PCM, the computation times for nested 
dissection and the Cholesky decomposition are negligible compared 
to the wavelet matrix compression. In view of the fill-in, there is 
an increase of approximately a factor of 5 for the average number 
of nonzero entries per row for \(N_J=1\,310\,720\). The computation 
times for the \(\theta\)-scheme in seconds for the 150 time steps are shown
in the column labelled by \(t_{\theta}\). As can be seen, we obtain a
solution time of roughly 11 minutes for \(N_J=1\,310\,720\) unknowns
in the spatial discretization.

\subsection{Gaussian random field}
\begin{figure}[htb]
\begin{center}
\includegraphics[scale=0.14,clip,trim= 360 75 360 75]{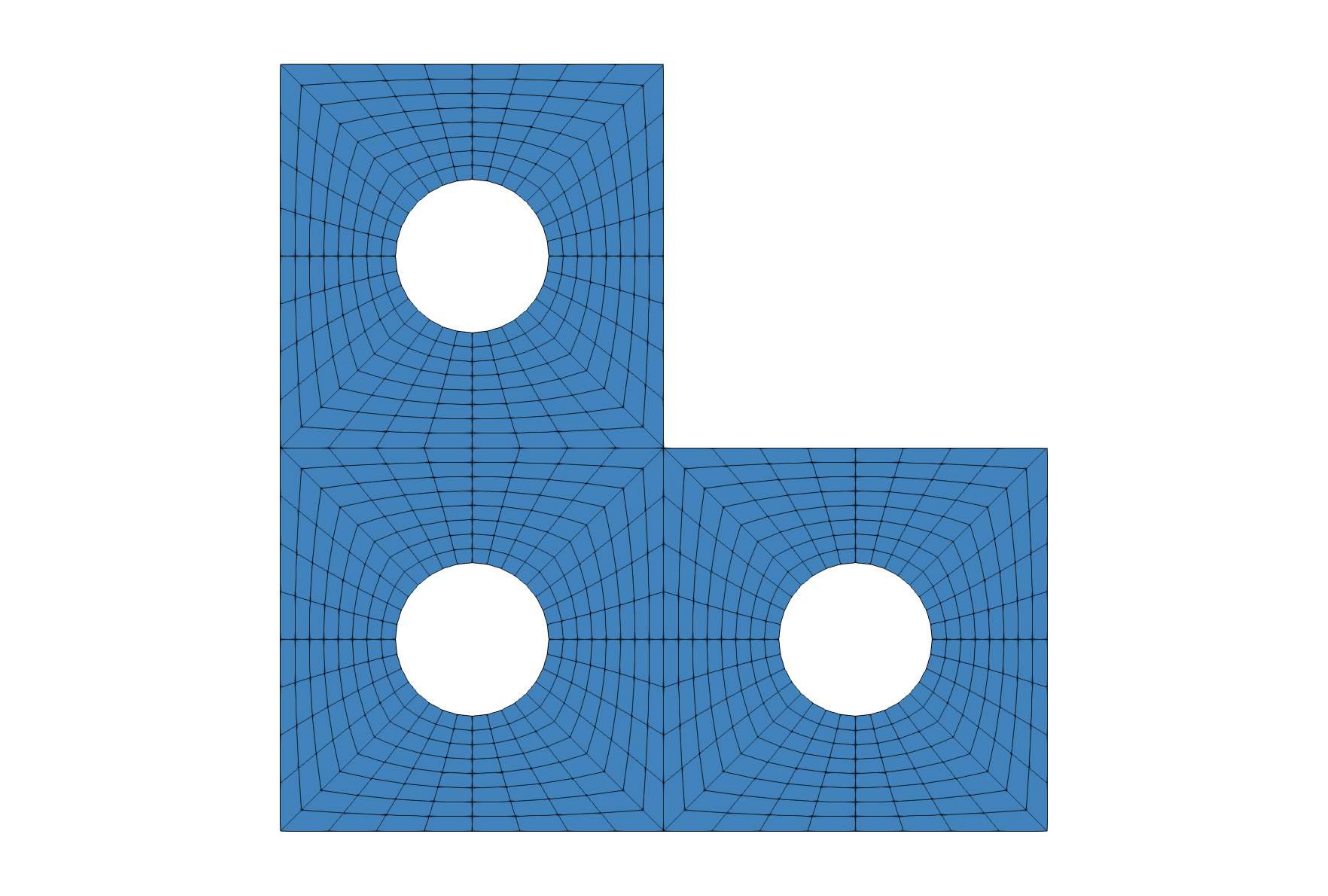}
\caption{\label{fig:randomFields}Computational geometry for the Gaussian random field.}
\end{center}
\end{figure}

\begin{figure}[htb]
\begin{center}
\begin{tikzpicture}
\draw(0,0) node {\includegraphics[scale=0.11,clip,trim= 380 120 610 60]{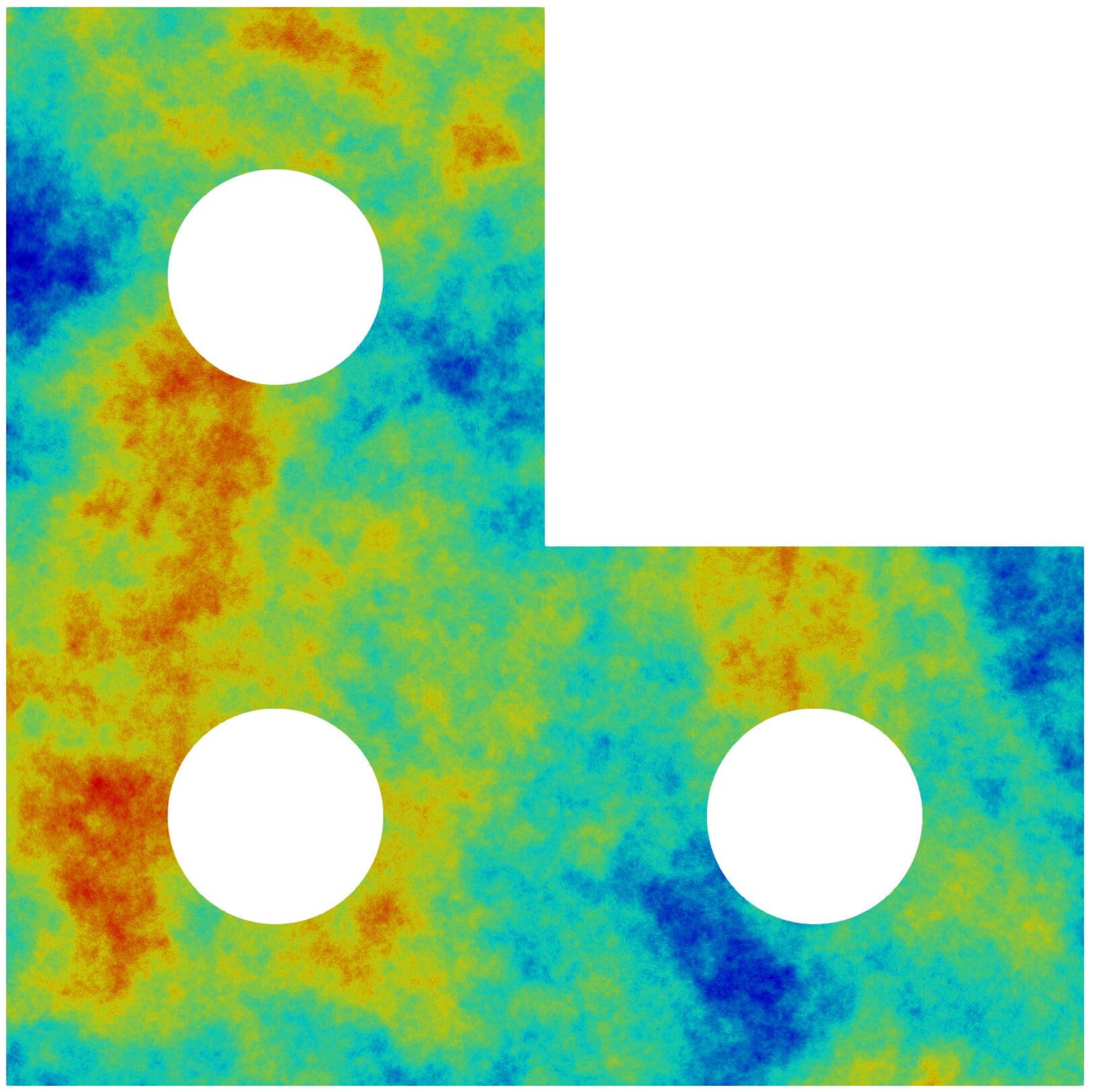}};
\draw(5,0) node {\includegraphics[scale=0.11,clip,trim= 380 120 610 60]{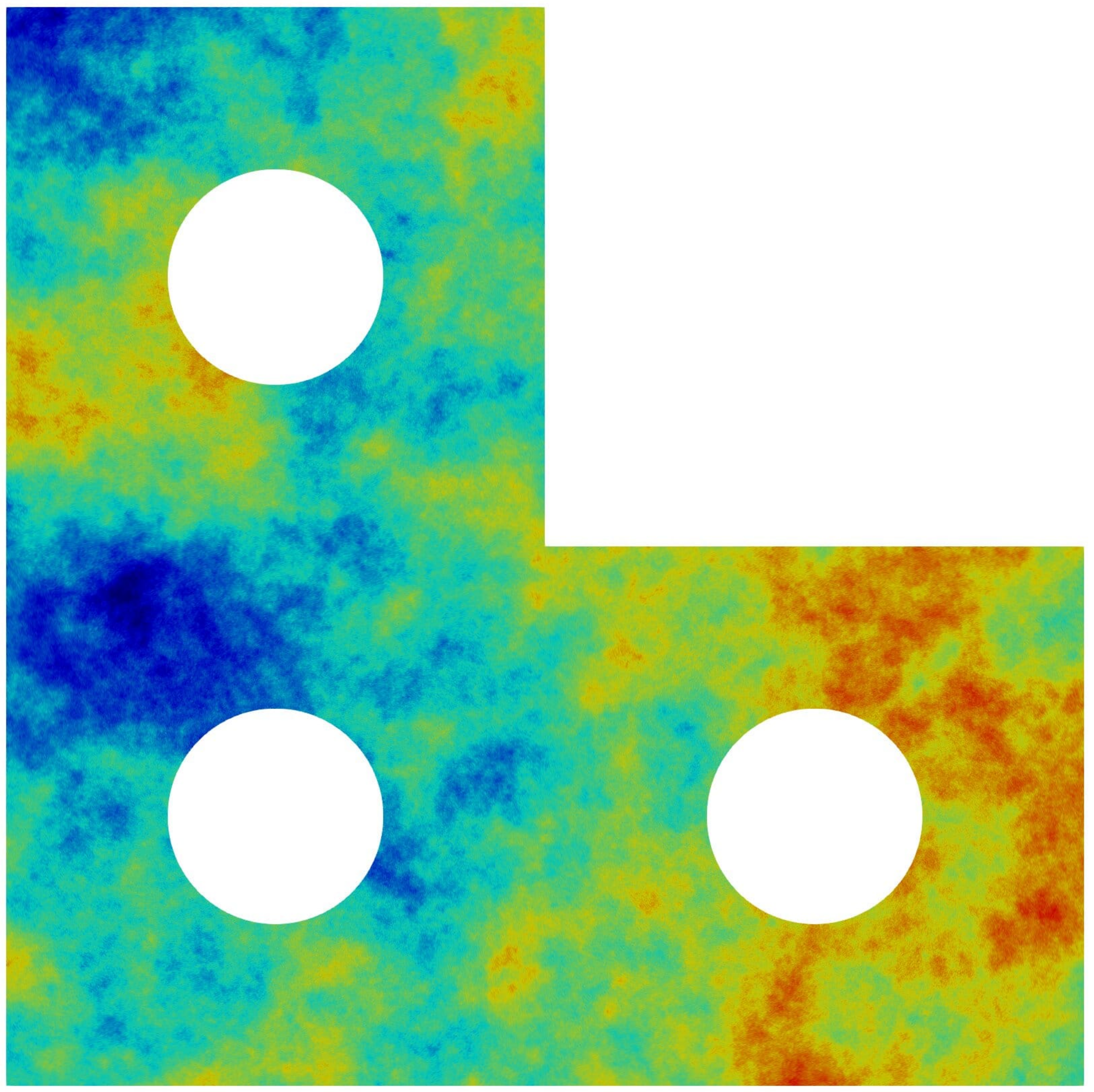}};
\draw(5,5) node {\includegraphics[scale=0.11,clip,trim= 380 120 610 60]{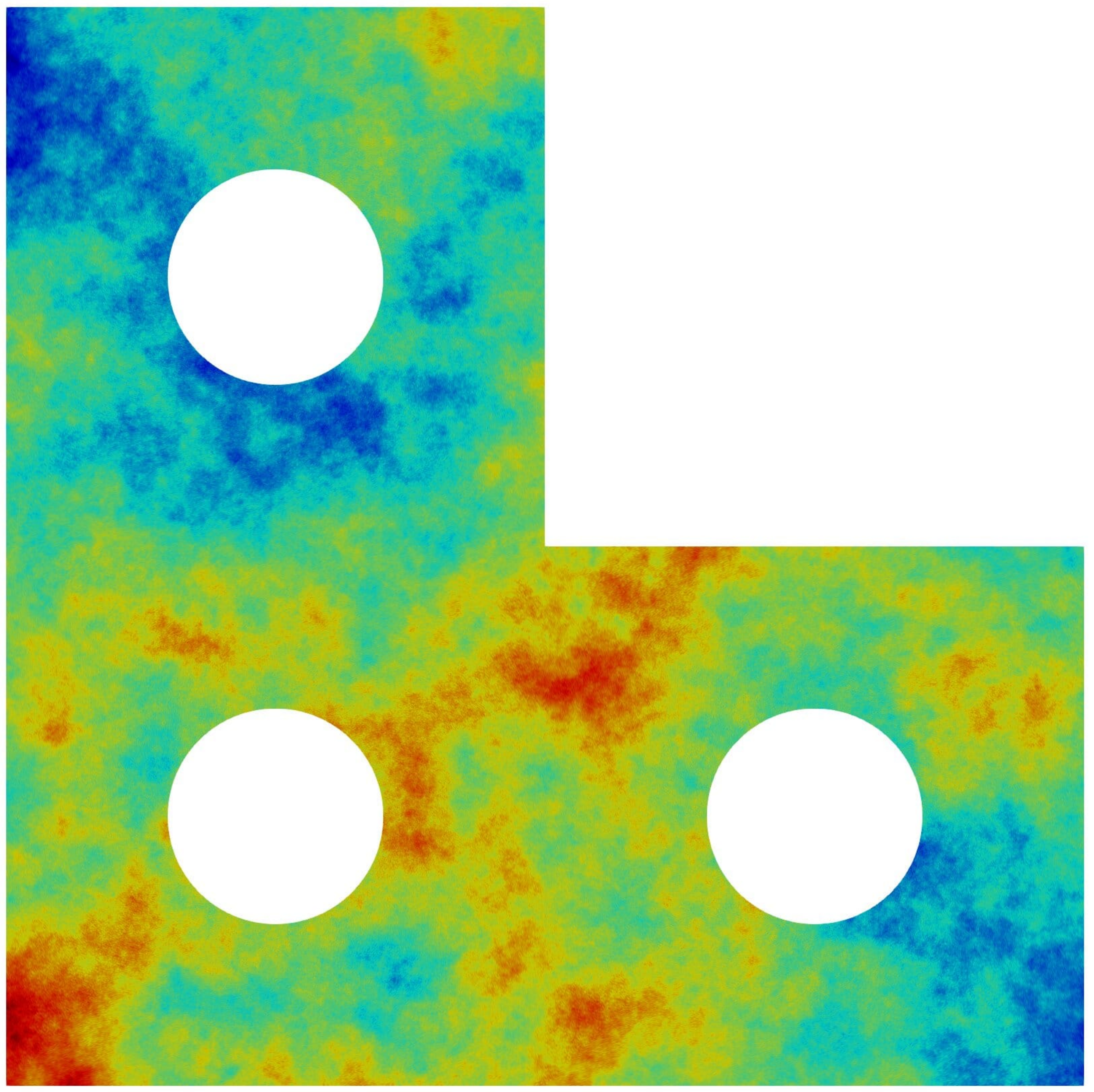}};
\draw(0,5) node {\includegraphics[scale=0.11,clip,trim= 380 120 610 60]{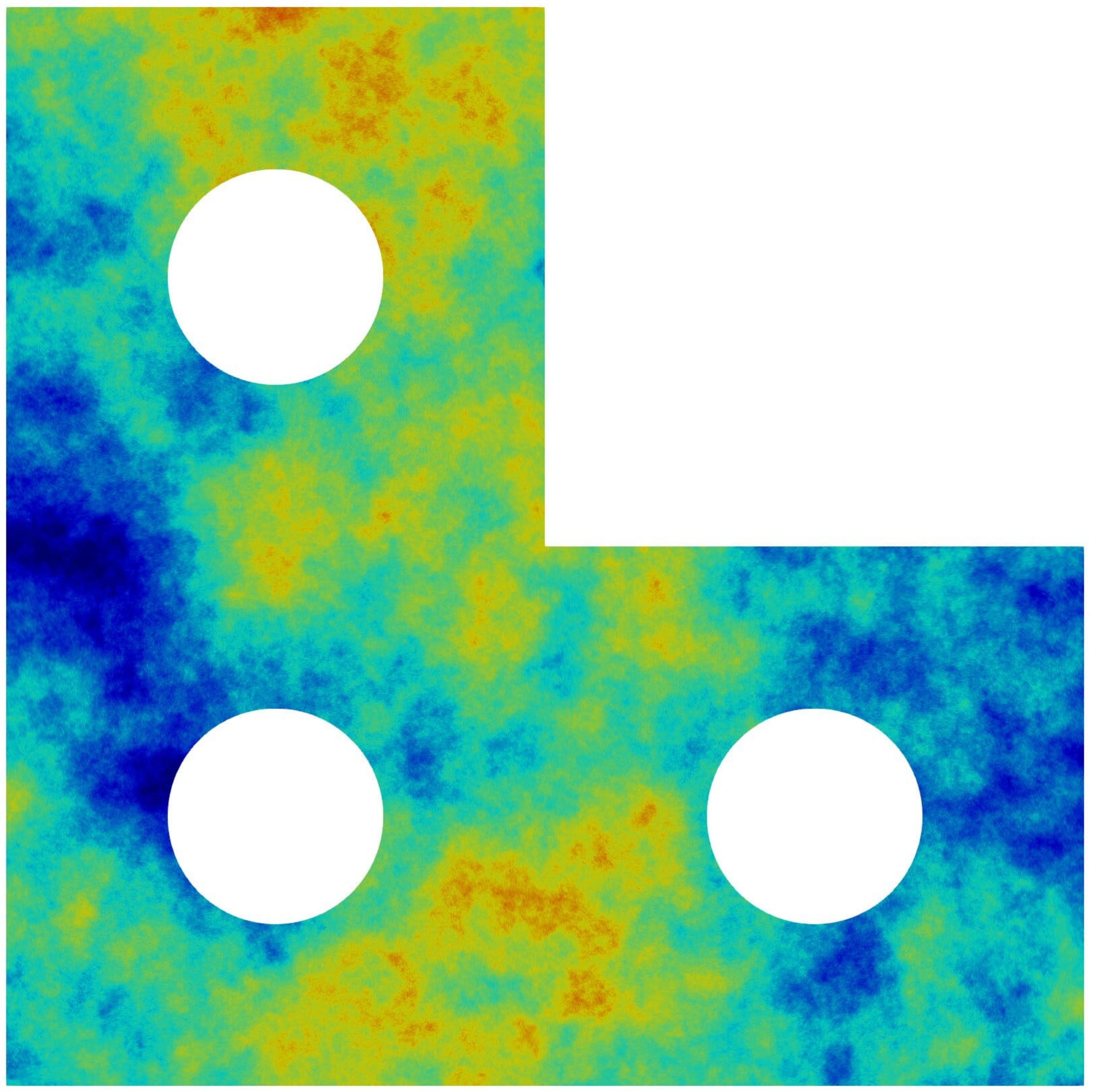}};
\draw(9.2,2.25) node {\includegraphics[scale=0.194,clip,trim= 1600 400 50 420]{field1}};
\fill[color=white] (8.5,0.5) rectangle (10.8,4);
\draw(9,3.8) node {$\phantom{-}3.8$};
\draw(9,2.98) node {$\phantom{-}2.0$};
\draw(9,2.54) node {$\phantom{-}1.0$};
\draw(9,2.1) node {$\phantom{-}0\phantom{.0}$};
\draw(9,1.65) node {$-1.0$};
\draw(9,1.22) node {$-2.0$};
\draw(9,0.75) node {$-3.0$};
\draw(10.1,2.25) node[rotate=90] {Value};
\end{tikzpicture}
\caption{\label{fig:FieldReals}Four different realizations of the Gaussian random field.}
\end{center}
\end{figure}

For the simulation of a Gaussian random field, we consider an
L-shaped domain with three holes, compare Figure~\ref{fig:randomFields}.
The domain has a side length of 4, while the holes have a diameter 
of 0.8. For the wavelet matrix compression, we use piecewise bilinear 
biorthogonal wavelets with four vanishing moments, see \cite{HS1}.
The expectation is set to \(\mathbb{E}[a]({\bs x})\equiv 0\), while the
covariance is given by the exponential kernel
\[
\operatorname{Cov}[a]({\bs x},{\bs x}')=\exp(-\|{\bs x}-{\bs x}'\|_2).
\]
Four different realizations of the corresponding field in case of 
\(N_J=792\,688\) nodes are depicted in Figure~\ref{fig:FieldReals}.

\begin{figure}[htb]
\begin{center}
\setlength{\fboxrule}{0.2pt}
\begin{minipage}{0.49\textwidth}
\begin{center}
\includegraphics[scale=0.55,clip,trim=0 0 0 13,frame]{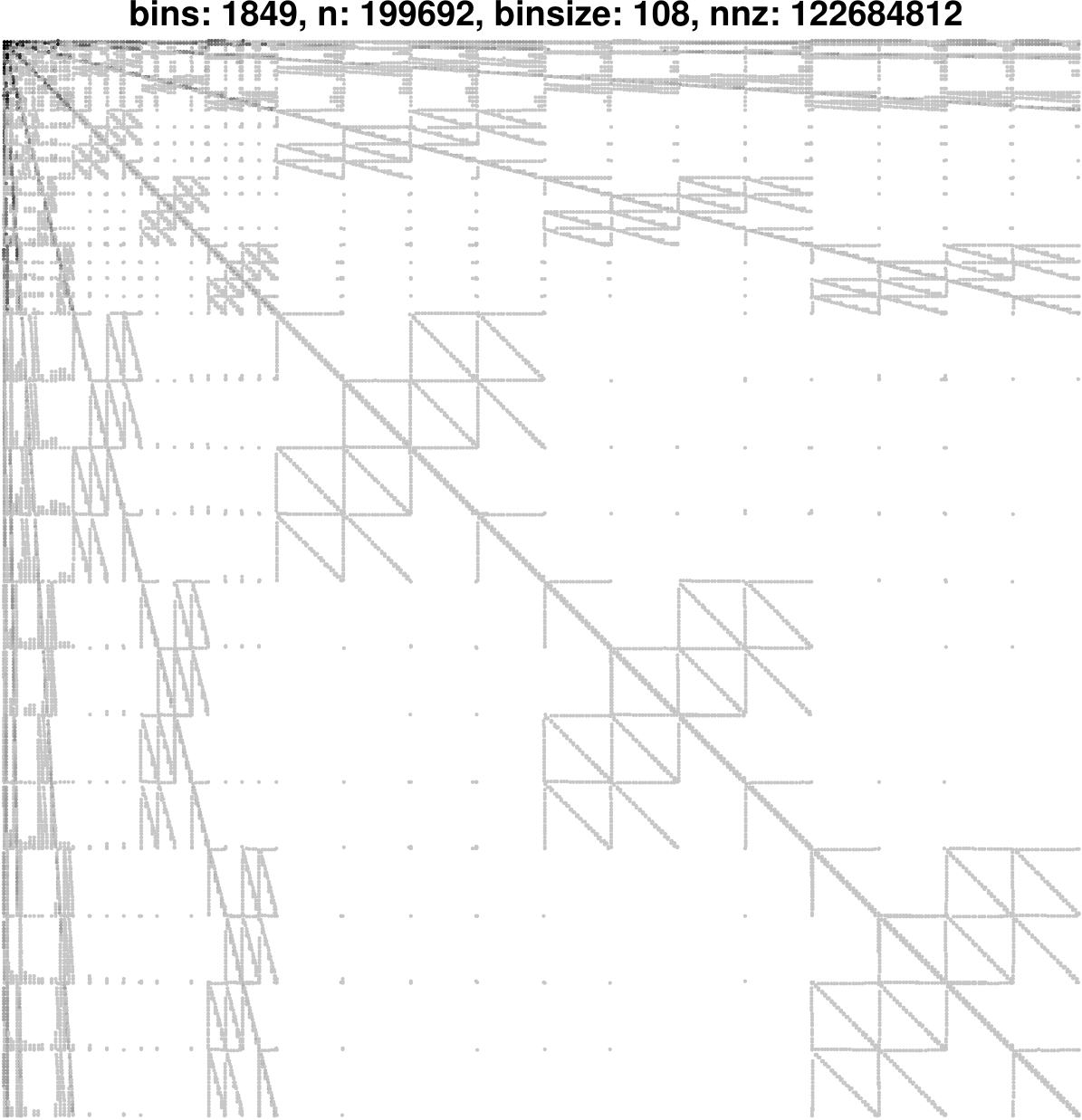}
$\nnz({\bs C}_J^G) = 122\,684\,812$
\end{center}
\end{minipage}
\caption{\label{fig:RandomFieldC}Sparsity patterns of \({\bs C}_J^G\) 
of the random field for \(N_J=199\,692\). Each dot corresponds to a submatrix
of size \(108\times 108\). Lighter blocks have less entries than darker blocks.}
\end{center}
\end{figure}

\begin{figure}[htb]
\begin{center}
\setlength{\fboxrule}{0.2pt}
\begin{minipage}{0.49\textwidth}
\begin{center}
\includegraphics[scale=0.55,clip,trim=0 0 0 13,frame]{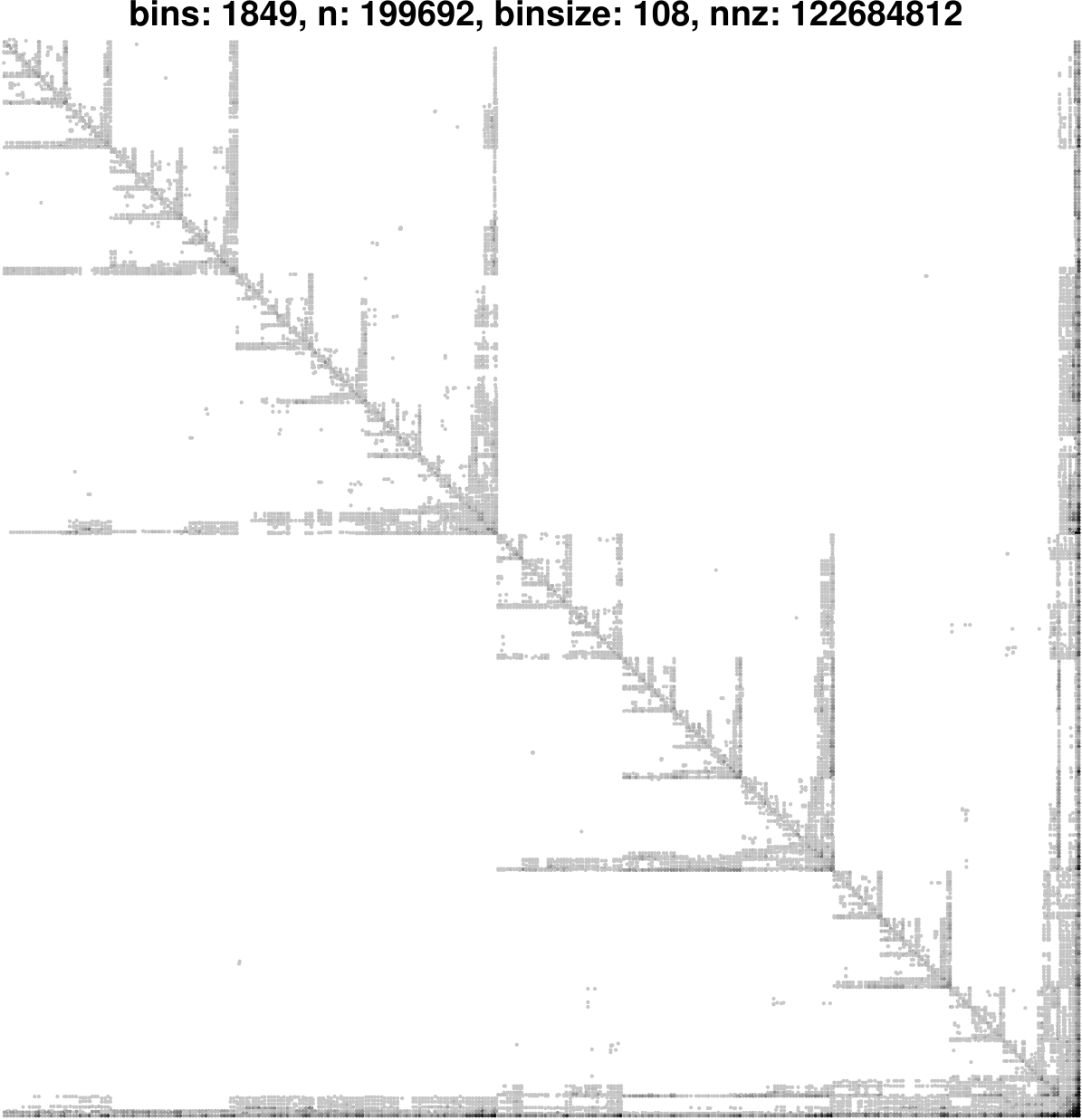}
$\nnz({\bs C}_J^G) = 122\,684\,812$
\end{center}
\end{minipage}
\begin{minipage}{0.49\textwidth}
\begin{center}
\includegraphics[scale=0.55,clip,trim=0 0 0 13,frame]{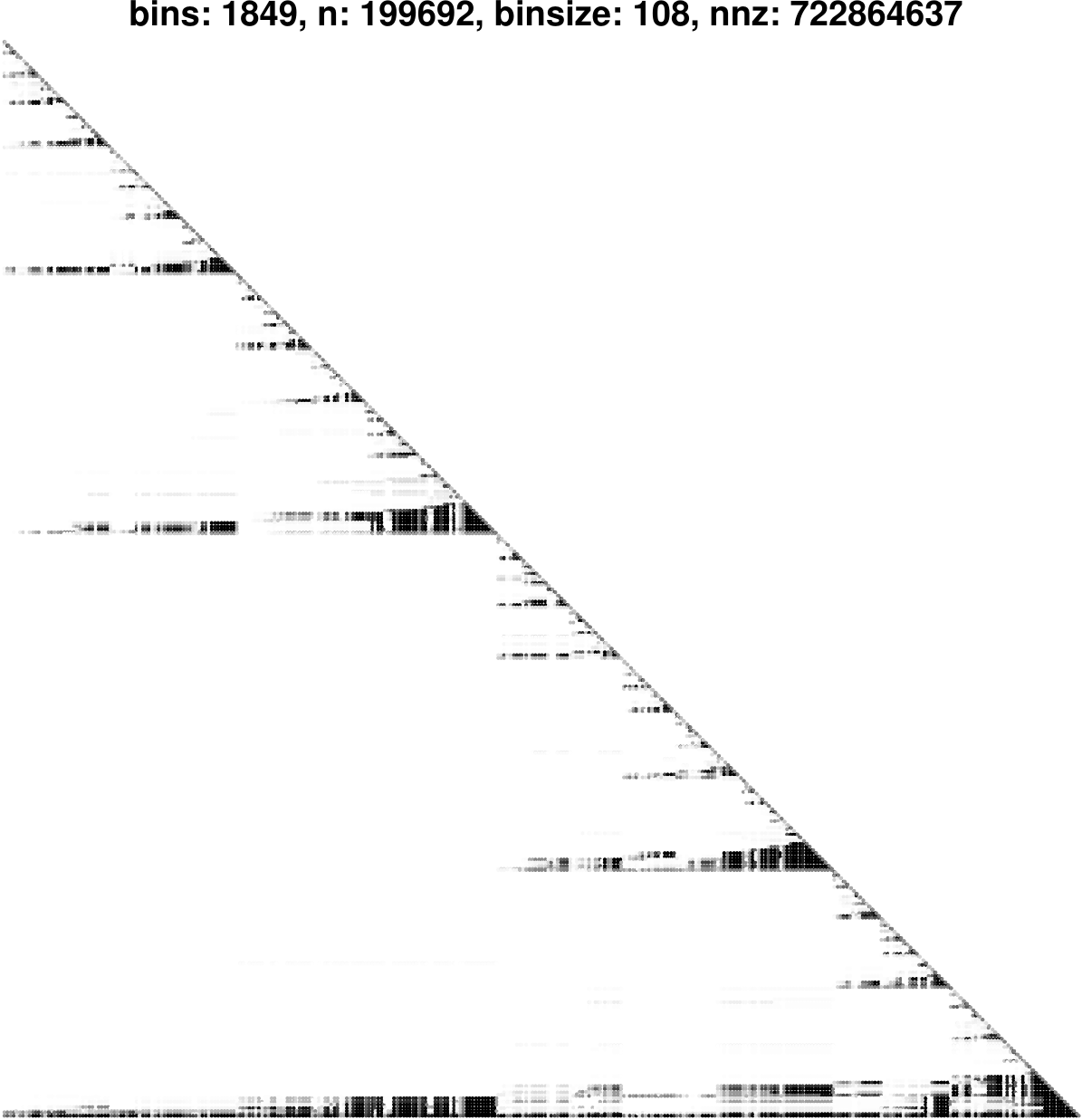}
$\nnz({\bs L}_J) = 722\,864\,637$
\end{center}
\end{minipage}
\caption{\label{fig:CholRF}Sparsity patterns of \({\bs C}_J^G\) (left) and 
the Cholesky factor \({\bs L}_J\) (right) for the random field and \(N_J=199\,692\). 
Each dot corresponds to a submatrix of size \(108\times 108\). Lighter blocks 
have less entries than darker blocks.}
\end{center}
\end{figure}

The pattern of the covariance matrix \({\bs C}_J^G\) is shown in 
Figure~\ref{fig:RandomFieldC}, while the patterns of the reordered 
matrix and its Cholesky factor are provided in Figure~\ref{fig:CholRF}.

Table~\ref{tab:RF} shows the computation times and the average 
numbers of nonzero entries as in the previous examples. In addition,
we have the column \(t_{\text{Sample}}\), which contains the times for 
computing a single realization of the random field in seconds. These 
times have been computed by averaging the computation times over 1000 samples.
In order to compute the inverse of the mass matrix, we could in
principle reuse the nested dissection ordering which has been computed for
the system matrix, as the pattern of the mass matrix is a subset of the 
pattern of the system matrix. However, this will result in a fill-in similar 
to the system matrix. Hence, it is favourable to use a different ordering 
for the mass matrix, resulting in much less fill-in.

\begin{table}[htb]
\begin{center}
\scalebox{0.93}{
\begin{tabular}{|r|r|r|r|r|r|r|r|r|}\hline
$N_J$ & $t_\text{WEM}$ &  $t_\text{ND}$ &   $t_\text{Chol}({\bs C}_J^G)$ &
$t_\text{Sample}$ & $\anz({\bs C}_J^G)$  & $\anz({\bs L}_J)$ &$\anz({\bs G}_J)$  &
$\anz(\widehat{\bs L}_J)$\\\hline
    972 &  1.50 & 0.04 &  0.23 &$6.25\!\cdot\!10^{-4}$  &      249  &     251
&  53 & 31\\ 
   3468 & 14.49 & 0.16 &  0.26 &$3.97\!\cdot\!10^{-3}$  &      338  &     498
&  91 & 72\\
  13068 &  119.44 & 0.83 &  2.22 &$3.13\!\cdot\!10^{-2}$  &      429  &    1096
& 126 &124\\
  50700 & 777.69 & 4.45 & 21.45 &$2.18\!\cdot\!10^{-1} $  &     517   &   1998
&156  &196\\
 199692 &  4368.40 &24.34 & 222.72 &  $1.51\!\cdot\!10^{0\phantom{-}} $  &     615   &   3620
&184  &283\\
792588 &   29376.29 &   130.23 & 4290.29 & $1.61\!\cdot\!10^{1\phantom{-}}$&   714 &   8336 &   211 &   390\\
\hline
\end{tabular}}
\caption{\label{tab:RF}Computation times and numbers of 
nonzero entries in case of the Gaussian random field.}
\end{center}
\end{table}

As can be seen from Table~\ref{tab:RF}, the times for sampling 
the random field only increase moderately. We remark that, due to the 
larger supports of the bilinear wavelets and the higher precision of the 
discretization, the system matrix contains more entries per row on average.
This also leads to higher computation times for the matrix assembly.
However, as before the increase of nonzero entries in the Cholesky factor 
remains very moderate. In addition, we have provided the average number 
of nonzeros per row for the mass matrix in the column labelled by 
$\anz({\bs G}_J)$. The number of nonzeros for the corresponding 
Cholesky factor of the mass matrix is given in the column labelled 
by $\anz(\widehat{\bs L}_J)$.

\section{Conclusion}\label{sec:conclusion}
In this article, we have proposed a very efficient direct
solver for nonlocal operators. The pivotal idea is to
combine the wavelet matrix compression with the nested dissection
ordering. Thereby, the fill-in resulting from a Cholesky
decomposition or an LU decomposition is drastically reduced. 
This has been numerically investigated into depth for three 
relevant applications, namely the polarizable continuum model, 
a parabolic problem for the fractional Laplacian in integral 
form, and the fast simulation of Gaussian random fields.
In all three cases, wavelet matrix compression yields sparse 
system matrices while the fill-in in the matrix factorization 
remains very low thanks to the nested dissection ordering. 
This behaviour could be observed for more then \(10^6\) 
unknowns in the discretization. 

A formidable application of the presented approach is the
fast simulation of rough (Gaussian) random fields. In such 
cases, the numerical solution of the eigenvalue problem for 
the covariance is computationally prohibitive and, hence, the 
computation of a Karhunen-Lo\`eve expansion is not feasible. 
In turn, the use of a wavelet basis yields a sparse representation 
of the covariance operator and a matrix root is rapidly computable 
by employing the nested dissection ordering and the Cholesky 
decomposition.

\subsection*{Acknowledgement}
This research has been funded in parts by the Swiss National Science 
Foundation (SNF) under project number $407540\_167186$ NFP 75 Big Data.
\bibliographystyle{plain}

\end{document}